\documentclass[11pt,letterpaper]{article}

\usepackage[T1]{fontenc}
\usepackage[utf8]{inputenc}
\usepackage[margin=1in]{geometry}
\usepackage{lmodern}
\usepackage{microtype}
\usepackage{setspace}
\usepackage{float}
\floatstyle{ruled}
\newfloat{algorithm}{htbp}{loa}
\floatname{algorithm}{Algorithm}
\usepackage{booktabs}
\usepackage{tabularx}
\usepackage{array}
\usepackage{amsmath,amssymb,amsthm,mathtools}
\usepackage{graphicx}
\usepackage{placeins}
\usepackage{needspace}
\usepackage[font=small,labelfont=bf,labelsep=period]{caption}
\usepackage{natbib}
\bibpunct[, ]{(}{)}{,}{a}{}{,}

\usepackage{hyperref}
\hypersetup{
  hidelinks,
  hypertexnames=false,
  pdftitle={Diagnosis and Recovery in AI Workflows with Scarce Expert Capacity},
  pdfauthor={Ruihan Zhou, Jiheng Zhang, Yijie Peng},
  pdfsubject={Sequential diagnosis and recovery under shared expert capacity},
  pdfkeywords={AI workflows; diagnostic information; checkpoint recovery; scarce expert capacity; operational control}
}

\newcolumntype{Y}{>{\raggedright\arraybackslash}X}
\newcolumntype{C}{>{\centering\arraybackslash}X}
\newcolumntype{L}[1]{>{\raggedright\arraybackslash}p{#1}}
\newcolumntype{P}[1]{>{\centering\arraybackslash}p{#1}}

\newcommand{\E}{\mathbb E}
\newcommand{\Prb}{\mathbb P}

\theoremstyle{plain}
\newtheorem{assumption}{Assumption}
\newtheorem{condition}{Condition}
\newtheorem{lemma}{Lemma}
\newtheorem{theorem}{Theorem}
\newtheorem{proposition}{Proposition}
\newtheorem{corollary}{Corollary}

\begin{document}

\begin{center}
{\LARGE\bfseries Diagnosis and Recovery in AI Workflows\\[3pt]
with Scarce Expert Capacity\par}
\vspace{12pt}
{\large Ruihan Zhou\textsuperscript{1}\qquad
Jiheng Zhang\textsuperscript{2}\qquad
Yijie Peng\textsuperscript{3}\par}
\vspace{8pt}
{\small
\textsuperscript{1}Guanghua School of Management, Peking University, Beijing 100871, China\\
\href{mailto:rhzhou@stu.pku.edu.cn}{rhzhou@stu.pku.edu.cn}\\[4pt]
\textsuperscript{2}Department of Industrial Engineering and Decision Analytics,\\
The Hong Kong University of Science and Technology, Hong Kong, China\\
\href{mailto:jiheng@ust.hk}{jiheng@ust.hk}\\[4pt]
\textsuperscript{3}School of Management and Engineering, Nanjing University, Nanjing, China\\
\href{mailto:pengyijie@nju.edu.cn}{pengyijie@nju.edu.cn}\par}
\end{center}
\vspace{5pt}

\begin{abstract}
\textbf{Problem definition:} We study sequential decisions in automated execution systems where failures of uncertain origin require recovery under limited shared expert capacity. Recovery choices affect future diagnostic demand, while expert availability changes which recovery actions are worthwhile. AI workflows provide a timely and representative setting for this operational problem. The objective is to coordinate diagnostic requests, expert dispatch, and checkpoint-based recovery over time to maximize expected net operating value.
\textbf{Methodology/results:} A finite-horizon model tracks uncertain error locations, remaining time, and shared expert capacity. Structural analysis shows how competition for diagnosis changes unaided recovery, without requiring an AI-specific execution technology. When a full rebuild from checkpoint 0 is the best unaided fallback after failure, a larger expert pool can justify retaining more work despite a lower chance of immediate success. Capacity can also raise net operating value before fixed staffing costs while reducing diagnostic hours when advice is requested only after a failed trial. Our Diagnosis--Resource--Recovery Control (DRRC) policy coordinates these decisions through local recovery lookahead and estimates of service requirements and remaining time; its learning-augmented variant uses task-specific quality and duration estimates. In controlled simulations, larger pools draw more requests before repair, and selective diagnosis raises net value chiefly by reducing recovery costs. In a digital twin based on 177 software defects, DRRC with reference-assisted diagnosis increases mean timely completion by 3.27 percentage points and mean net value by 1.85\% relative to automated recovery, although the net-value advantage remains uncertain across repositories.
\textbf{Managerial implications:} Expert capacity changes when workflows seek advice and how much work they preserve. Expansion may raise or lower diagnostic hours, so capacity and recovery design should be assessed jointly by their contribution to net operating value.
\end{abstract}


\section{Introduction}
\label{sec:introduction}

AI agents are extending automation from individual responses to connected stages of work, changing how organizations use human expertise. As systems plan, use tools, and revise outputs, specialists direct execution, assess uncertain results, and resolve exceptions. Evidence from Anthropic documents longer autonomous sequences of coding work alongside continued demands on human oversight \citep{huang2025aiwork}. The ability to execute more tasks does not by itself ensure that more valuable work is completed on time. Failures can require repeated execution or diagnosis, and several workflows may need the same specialist. Organizations must therefore coordinate automated work and scarce expertise across tasks, accounting for the additional work generated by their own recovery decisions.

Current platforms provide practical building blocks for this division of work. GitHub Copilot can plan changes, edit code, and run tests \citep{github2026copilot}, while AWS Step Functions can resume unsuccessful executions without repeating successful steps \citep{aws2026redrive}. These capabilities support automated execution and recovery that preserves completed work. AI workflows are a particularly timely and representative setting for this operational problem: execution spans several stages, verification may leave the source of an error uncertain, and recovery can use both automated feedback and specialist advice. For example, in a software-maintenance workflow, an agent interprets an issue, edits code, and runs tests. A failed test may originate in the requirements, the patch, or a dependency. Retrying the last step may preserve an earlier error, whereas restarting the entire task discards useful work. Expert diagnosis can help select a suitable restart point and repair procedure, but it takes time and delays other requests. Before submitting a request, the controller can instead let the agent try a repair while a specialist is busy. Success saves a diagnostic visit. If the repair fails, the controller gains evidence about the fault, but other workflows may then need the same specialist. The prospect of obtaining advice later therefore matters even to a workflow that proceeds on its own today.

This interaction makes diagnostic demand endogenous to the operating policy. Incoming workflows create potential needs for assistance, but whether and when they enter the expert queue depends on recovery choices and the evidence those choices produce. Requests change congestion; congestion changes the time and opportunity cost of obtaining advice; and access to advice changes which work is retained during recovery. The resulting outcomes then determine subsequent requests. A larger expert pool can therefore change the demand it serves even when incoming workflows and recovery technology remain unchanged. The fundamental problem is how to operate an automated execution system with this dynamic feedback so as to maximize net operating value, accounting for timely completion, expert service, failures, delay, and repeated work.

\begin{figure}[htbp]
\centering
\includegraphics[width=1\textwidth]{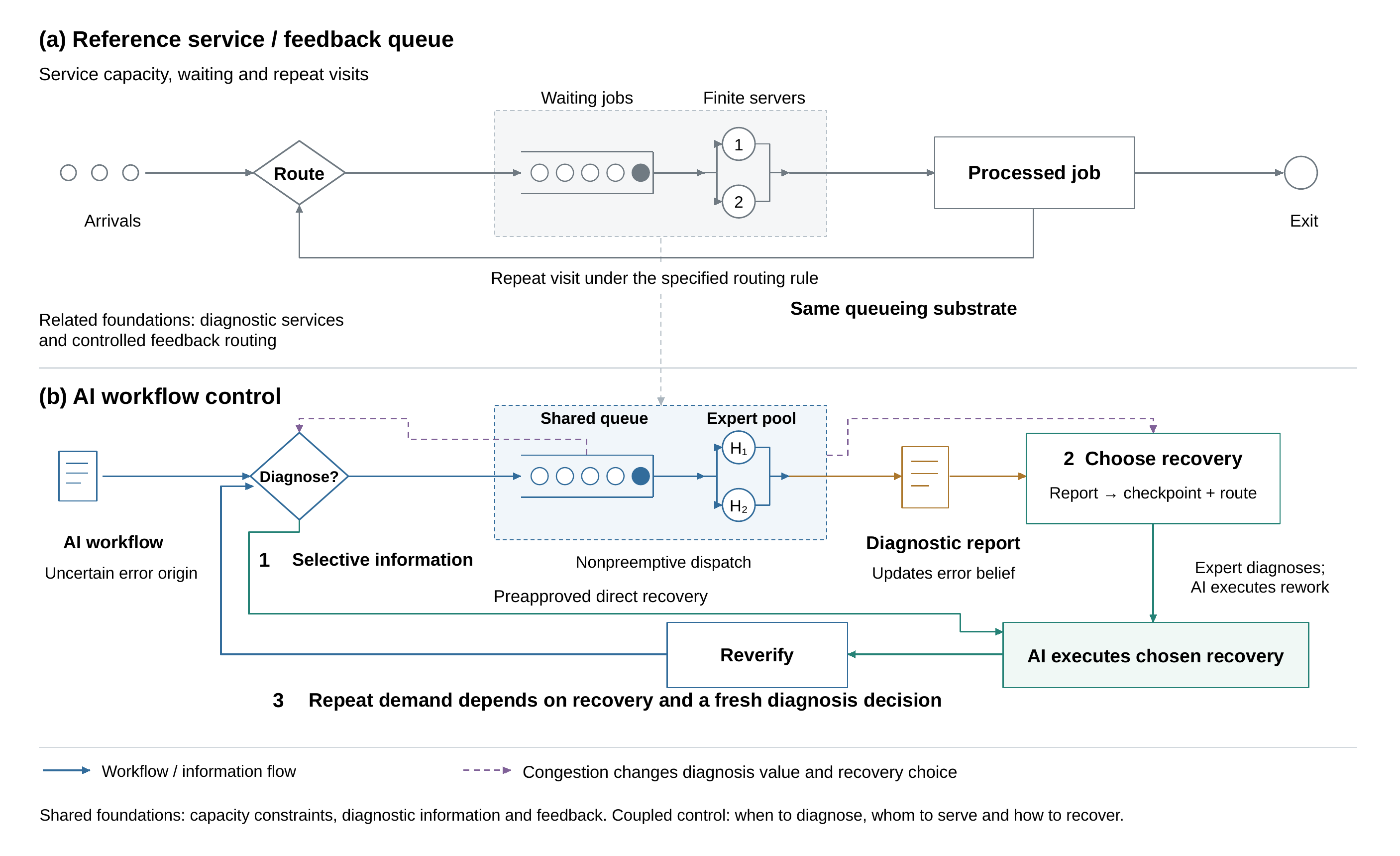}
\caption{Shared queueing structure and coordinated diagnosis--recovery decisions.}
\label{fig:framework}
\end{figure}

Effective use of expertise requires deciding when to acquire information, how to allocate the capacity that supplies it, and how to act on it. We therefore coordinate three linked decisions: whether a workflow should request diagnosis at screening (E2), which waiting request a free expert should serve (E3), and which recovery plan to execute after a report (E4). Automated verification supplies inexpensive but imperfect evidence. The expert adds diagnostic information, the controller selects a saved checkpoint and recovery route, and the execution agent carries out that plan. The repaired output returns to verification, and another expert visit requires a fresh request after subsequent screening. Any required formal certification remains a separate release gate. Coordinating these decisions determines whether advice arrives early enough to improve recovery and whether the value of the work saved justifies the resources used. Figure~\ref{fig:framework} places this interaction alongside a reference service and feedback queue. Both feature limited servers, waiting, and repeat visits; here a report informs a checkpoint--route choice whose outcome affects subsequent demand. The same operational structure can arise in cloud service recovery and automated troubleshooting when completed work can be preserved or repeated and shared specialists supply information for recovery. AI workflows thus provide a concrete setting for a broader operations-management problem. This connects selective human intervention under congestion \citep{lykouris2025defer} with repeat-service control \citep{chan2024returns} through information-dependent recovery.

Our first contribution is to show how competition for diagnosis changes the work a platform preserves and the timing of expert intervention. We model the three decisions in Figure~\ref{fig:framework} as a finite-horizon belief-state semi-Markov decision process (SMDP). To identify how access to later advice changes the initial recovery choice, we then analyze a simpler system that isolates competition after a recovery attempt. A specialist may be available now but face other requests when that attempt ends. Adding an expert can lower the opportunity cost of diagnosis after failure. When a full rebuild from checkpoint 0 is the best unaided fallback, the controller can retain more completed stages despite a greater risk of another repair. It can also switch from immediate diagnosis to seeking advice only on failure, raising net operating value before fixed staffing costs while reducing expert time. The structural results apply to automated execution systems satisfying the stated information, recovery, and capacity conditions, whether execution is AI-based or conventional.

Our second contribution is a controller that coordinates immediate and contingent diagnosis with recovery choices in the general system. Exact SMDP optimization must track test results, reports, and service completions across concurrent workflows to evaluate a recovery decision. Diagnosis--Resource--Recovery Control (DRRC) instead looks ahead within each workflow and represents competing demand through predicted waiting and a charge per unit of expert time. Screening compares immediate diagnosis with acting on current evidence, including the option to request advice after another test. Dispatch accounts for the time left to use a report, and recovery uses the information actually obtained. LA-DRRC uses historical records to estimate task-specific recovery quality and service time within the same control structure.

Controlled experiments show that LA-DRRC earns higher net value than optimized automated recovery and one-step diagnostic control, chiefly by reducing recovery costs, and higher average value than tree search (UCT) and reinforcement learning (PPO) on independent test episodes. Varying checkpoints and report accuracy identifies when diagnosis saves useful work; disabling selective requests or report use reduces policy value. A digital twin combining observed repair outcomes for 177 defects from 11 repositories with a simulated queue shows higher average timely completion under DRRC with reference-assisted diagnosis than under automated recovery. The application evaluates selective access to advice and report-dependent recovery under shared capacity; the controlled experiments also examine checkpoint choices and repeated diagnostic visits. They also show why capacity planning must consider changes in recovery behavior: historical request volumes reflect choices made under the existing expert pool and need not predict demand after expansion.

Section~\ref{sec:literature} positions the study in the related literature. Section~\ref{sec:model} specifies the operational model and exact SMDP. Section~\ref{sec:capacity} isolates the effect of competition on recovery, and Section~\ref{sec:drrc} uses this analysis to construct an approximate policy for the general model. Section~\ref{sec:experiment} presents the synthetic evaluation and empirical application, and Section~\ref{sec:conclusion} discusses the managerial implications and future research.

\section{Related Literature and Research Positioning}
\label{sec:literature}

Our study connects three strands of research: AI--human intervention under congestion, diagnostic information and repeat service, and prediction for operational control. They meet at the choice between recovery on current evidence and requesting shared diagnostic expertise.

\subsection{AI workflows and human intervention}

Research on AI at work documents productivity gains and heterogeneity across workers and tasks \citep{cui2026generative}. Agent systems combine reasoning, tool use, and feedback across attempts \citep{yao2023react,shinn2023reflexion}, with repository-level systems and benchmarks extending these capabilities to software engineering \citep{jimenez2024swebench,yang2024sweagent}. Checkpoint/restore runtimes and rollback-based correction support recoverable execution \citep{wu2026crab,hao2026rollback}. We take execution and recovery capabilities as primitives and study how to coordinate expert attention when several workflows encounter uncertain failures.

Models of AI--human service systems connect intervention to congestion. \citet{lykouris2025defer} jointly study classification, admission to human review, and scheduling, with review correcting errors and supplying delayed, selectively observed feedback. \citet{lin2026screen} relate optional AI screening before mandatory approval to resource bottlenecks and give a fluid-limit guarantee; their extension allows human feedback to improve revisions. AI workflow queue models show how review and rework can absorb gains from faster execution \citep{bartolucci2026queue}. Here, expert intervention supplies information for a subsequent checkpoint--route choice; mandatory approval remains separate. This sequence connects human attention to the work retained during recovery and the diagnostic demand generated by later failures, requiring joint treatment of information value and repeat service.

\subsection{Diagnostic information, rollback, and repeat service}

Diagnostic-service research relates information acquisition to capacity: \citet{alizamir2013diagnostic} characterize how congestion changes the optimal stopping rule for sequential testing. In decision-theoretic troubleshooting, an attempted repair supplies evidence for the next test or intervention \citep{breese1996troubleshooting}; rollback recovery trades re-execution cost against the risk of retaining damaged work \citep{lin1998damage}. With a shared specialist, the opportunity to investigate an unsuccessful rollback depends on competing demand, linking information acquired through repair to the choice of how far to roll back.

A report's value depends on the available recovery actions, consistent with the decision-flexibility analysis of \citet{merkhofer1977flexibility}. Our controlled comparisons vary checkpoints and report accuracy separately to distinguish cheaper direct recovery from better use of advice. Feedback queues describe repeat service \citep{davignon1977feedback}, while service-speed and postservice-intervention models show how current actions affect returns \citep{huang2015patientflow,chan2014speedup}. \citet{chan2024returns} balance intervention costs against later service and congestion costs. We apply this lifecycle perspective to information-dependent recovery: reports change checkpoints and routes, recovery changes failure risk, and subsequent screening determines new requests. Expert demand thus follows from recovery and information rather than a directly selected return probability. Scheduling with abandonment and patient-flow models with deadlines and feedback provide foundations for dispatch under delay \citep{long2020scheduling,huang2015patientflow}. Here, waiting can leave too little time to use a report, making dispatch and recovery decisions interdependent.

\subsection{Prediction and operational control}

Decision-focused learning evaluates predictions through optimization outcomes \citep{elmachtoub2022smart}, while learning-augmented algorithms combine predictive inputs with explicit decision rules \citep{lykouris2018caching}. LA-DRRC uses supervised estimates within a prescribed controller: recovery-quality estimates distinguish plans, and diagnostic-time estimates determine whether advice can arrive in time to be useful. Other queueing methods learn decisions or system primitives through differentiable simulation \citep{che2024differentiable}, joint pricing and capacity learning \citep{chen2023pricing}, or online optimization with unknown arrival and service parameters \citep{chen2026unknown}. Our focus is on predicting recovery outcomes and diagnostic duration for an explicit control rule.

\section{Operational Model}
\label{sec:model}

We model a platform allocating shared experts across multistage workflows with deadlines. An execution agent performs and recovers work from saved checkpoints, an automated verifier supplies imperfect error evidence, and an expert provides diagnostic information. The agent may be an AI system or a conventional automated program. The model represents systems in which completed work can be preserved or repeated and recovery outcomes can generate further requests for shared diagnostic resources. AI workflows provide the motivating setting. Workflow arrivals are specified externally; diagnostic requests arise from the controller's decisions and the outcomes of execution and recovery. The controller observes the available evidence and the queue, balancing timely completion against expert service and repeated-work costs.

\subsection{Workflow dynamics and decision epochs}
\label{sec:modelarchitecture}

\paragraph{Workflow structure.}
Workflows arrive asynchronously over a finite horizon $[0,T]$. Let $\mathcal I(t)$ be the workflows currently inside the system at time $t$. Workflow $i$ has deadline $D_i$, completion value $v_i$, escaped-error loss $\ell_i^E$, and deadline or abandonment loss $\ell_i^D$. It contains stages $1,\ldots,J_i$. A stage creates a persistent work product, such as a plan, selected data, a model file, a code revision, a test result, or a technical document.

The workflow platform saves some intermediate states. Let $\mathcal K_i(t)$ be the set of saved states that can currently be restored. Checkpoint $k\in\mathcal K_i(t)$ stores the artifact, execution context, and dependencies after stage $k$ and before stage $k+1$. Checkpoint availability is exogenous: the controller chooses among saved states supplied by the platform. Let $c_i$ be the last checkpoint that passed formal certification and let $j_i$ be the current decision gate. Then $c_i<j_i$ and $\mathcal K_i(t)\subseteq\{c_i,\ldots,j_i-1\}$.

A recovery plan $u=(k,g)$ selects a saved checkpoint $k$ and a route $g\in\mathcal G_i(k)$. The checkpoint determines how much completed work is retained; the route specifies how the discarded work is rebuilt and validated. A rebuild route from $k>0$ still retains the prefix through $k$; a full rebuild starts from checkpoint 0. Depending on the application, it may change the repair instruction, model or agent, tools, or validation protocol. The organization supplies the feasible route set at each checkpoint, and the controller selects among these alternatives.

\begin{figure}[htbp]
\centering
\includegraphics[width=0.95\textwidth]{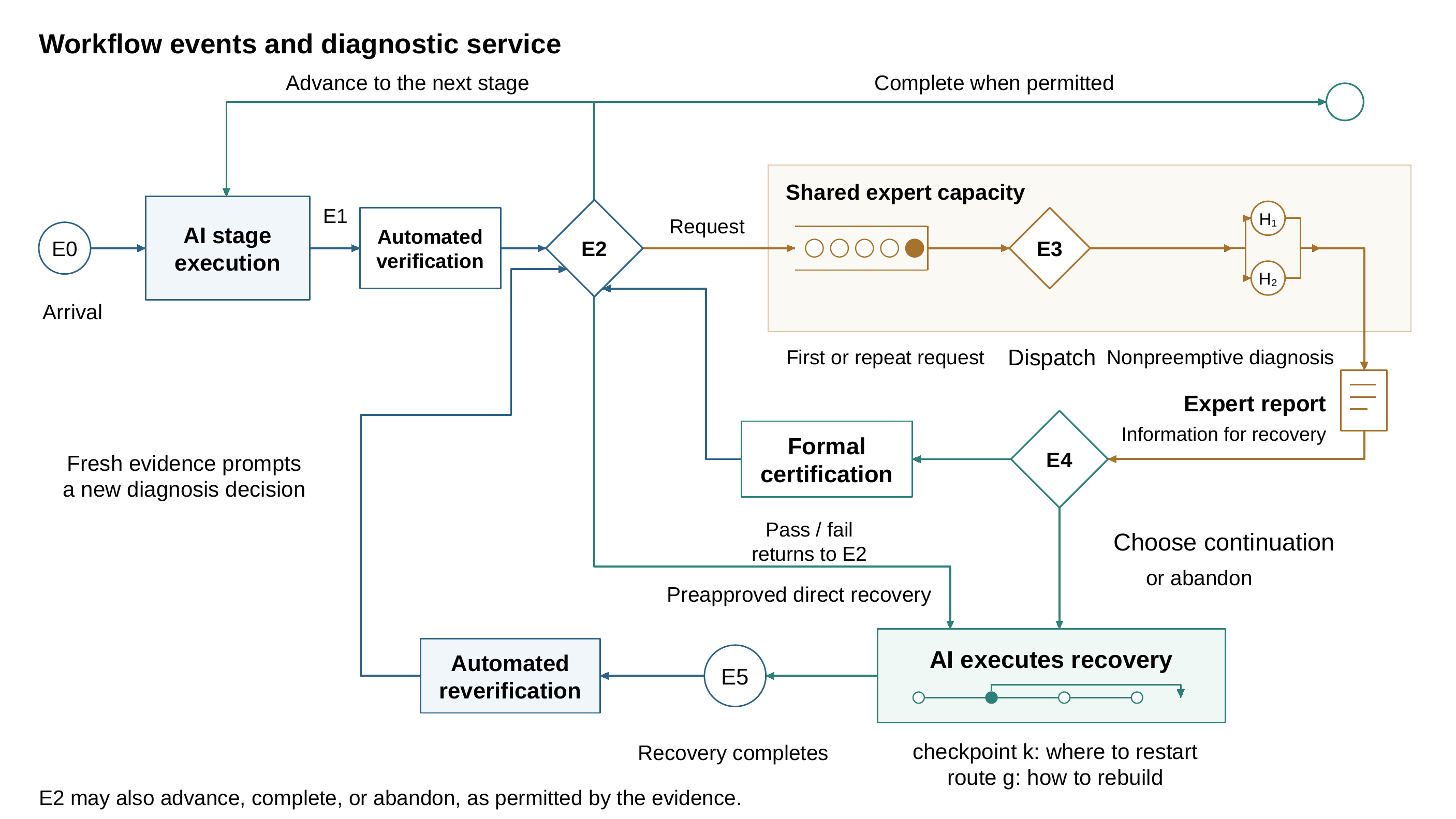}
\caption{Workflow events and decision points (E0--E5). E2 also permits direct submission to certification when its prerequisites are satisfied.}
\label{fig:workflow}
\end{figure}

For illustration, consider a three-stage workflow for data preparation, analysis, and recommendation: $u=(2,g_{\mathrm{target}})$ retains data and analysis and reruns recommendation generation; $u=(1,g_{\mathrm{rebuild}})$ retains only the data and rebuilds the later stages. Figure~\ref{fig:workflow} shows events E0--E5; these event labels are distinct from stage indices.

\paragraph{Arrival, execution, and screening (E0--E2).}
At E0 the system admits the workflow and records its attributes, value, and deadline. At E1 a completed stage goes to verification. The evidence at E2 determines whether the controller may advance, complete, submit the output to certification, recover directly, request diagnosis, or abandon. A diagnostic request joins $\mathcal Q(t)$ and triggers E3 if an expert is free. Once submitted, it remains queued until service or expiry; recovery does not proceed while the request is queued.

\paragraph{Expert dispatch (E3).}
A new request finding an available expert or a service completion triggers E3. The controller assigns a feasible waiting request or idles. With several free experts it updates the queue after each assignment. Once started, service is nonpreemptive.

\paragraph{Diagnosis and recovery completion (E4--E5).}
At E4 a report makes it possible to select a checkpoint--route plan using the new evidence, submit the work to certification, or abandon. When the execution agent completes the selected repair at E5, the verifier tests the result and the controller faces a new E2 decision. A return to the expert queue requires a fresh request.

\subsection{Hidden error state and information updates}
\label{sec:modelbelief}
\begin{table}[htbp]
\centering
\caption{Core Notation\label{tab:corenotation}}
{\scriptsize\renewcommand{\arraystretch}{1.08}\begin{tabularx}{\textwidth}{P{0.16\textwidth}L{0.29\textwidth}P{0.16\textwidth}Y}
\toprule
Symbol & Meaning & Symbol & Meaning\\
\midrule
$j_i,c_i,\mathcal K_i(t)$ & Current gate, last certified checkpoint, and restorable checkpoints & $Z_i,\mathbf b_i(t)$ & Earliest unresolved error and its posterior belief\\
$\mathbf y_i,\mathbf d_i$ & Current verifier/certification record and expert report & $u=(k,g)$ & Recovery checkpoint and recovery route\\
$\xi_i,\tau_i$ & Operating phase and remaining deadline slack & $\delta_i^X,\delta_i^V,\delta_i^R$ & Observable execution, verification, and recovery timing\\
$\mathcal Q(t),B$ & Waiting diagnostic requests and number of experts & $S_i^H,H_i(\pi)$ & One diagnostic service time and lifecycle expert work\\
$\mathbf s_i(t),\mathbf S(t)$ & Workflow state and full physical system state & $\pi,V_n$ & Operating policy and exact continuation value\\
$\mathcal U_i^A,\mathcal U_i^H$ & Repair plans permitted before and after diagnosis & $\mathcal Q^F(t)$ & Requests physically eligible for service\\
$q_k,\mathbf b^k$ & Failure probability and posterior after a trial from $k$ & $\theta_B,\theta^\star$ & Value of the report displaced at capacity $B$ and the cutoff for trying a repair first\\
$\mu$ & DRRC charge per unit of expert time & $\Delta_i^H$ & DRRC net gain from requesting diagnosis now\\
\bottomrule
\end{tabularx}\par}

\end{table}
Table~\ref{tab:corenotation} summarizes the model notation and the principal quantities used in Sections~\ref{sec:capacity} and~\ref{sec:drrc}. The platform can save work after any eligible stage. It treats a prefix as trusted only after formal certification. Certification is a separate verification activity and does not consume the diagnostic capacity $B$.

\begin{assumption}[Certified prefixes]\label{ass:certification}
A formal certification pass establishes that the work through the certified checkpoint is valid. Routine checks and expert reports may be imperfect.
\end{assumption}

Uncertainty therefore concerns work downstream of the last certified checkpoint $c_i$. The hidden state is the stage of the earliest unresolved error,
\[
Z_i\in\mathcal Z_i(t)=\{0,c_i+1,\ldots,j_i\},
\]
where $Z_i=0$ means that the current artifact is valid. If $Z_i=z>0$, the
earliest unresolved error occurred at stage $z$. Given the workflow history
$\mathcal H_i(t)$, including the common operating conditions observed by time $t$, the controller maintains the posterior belief
 \begin{equation}
\mathbf b_i(t)
=
\bigl(b_{iz}(t)\bigr)_{z\in\mathcal Z_i(t)},
\qquad
b_{iz}(t)
=
\Prb\{Z_i=z\mid\mathcal H_i(t)\}.
\label{eq:belief}
\end{equation} 
The earliest unresolved error determines how far recovery must roll back; conditional kernels summarize subsequent errors under Assumption~\ref{ass:markov}. Error location matters beyond failure probability: a stage-1 error requires checkpoint 0, whereas a stage-3 error can be removed from checkpoint 2, potentially at a different cost.

\paragraph{Verification and expert information.}
The verifier record $\mathbf y_i=(y_i^{\mathrm{pass}},y_i^{\mathrm{score}},\mathbf y_i^{\mathrm{loc}},\mathbf y_i^{\mathrm{trace}})$ contains the test result, score, location signals, and trace evidence, and determines preapproved direct plans $\mathcal U_i^A(\mathbf y_i)$. The report $\mathbf d_i=(d_i^{\mathrm{valid}},\mathbf d_i^{\mathrm{loc}},d_i^{\mathrm{severity}},d_i^{\mathrm{confidence}},\mathcal R_i^{\mathrm{feasible}})$ records validity, likely error location, severity, confidence, and feasible repairs. Diagnosis can change both the belief and the authorized set $\mathcal U_i^H(\mathbf y_i,\mathbf d_i)$.

Route $g$ at stage $j$ has transition kernel
$P_i^{X}(z',j'\mid z,j,g)$. The automated verifier has signal kernel
$P_i^{V}(\mathbf y\mid z,j,g)$, where
$\mathbf y\in\mathcal Y_i$ is the structured evidence observed at E2. A
recovery plan $u=(k,g)$ has transition kernel
$P_i^{R}(z',j'\mid z,k,g)$, where $z'$ and $j'$ denote the next error state and decision gate. We impose
\begin{equation}
P_i^{R}(z',j'\mid z,k,g)=0
\quad
\text{if }1\leq z\leq k\text{ and }z'\neq z.
\label{eq:rollbackrestriction}
\end{equation}
Equation~\eqref{eq:rollbackrestriction} states that rollback to checkpoint $k$
cannot remove an unresolved error that occurred at or before stage $k$. To
remove an error at stage $z$, the recovery plan must therefore restart from a
checkpoint strictly before $z$. If $z=0$ or $z>k$, re-execution can produce a valid artifact or introduce an error in a re-executed stage, whose index is strictly greater than $k$.

Execution or recovery can change both the error location and the next gate. Let $P_i^a(z',j'\mid z,j)$ be their joint transition law for action $a$, using the route and, for a repair, the checkpoint selected by $a$. For an observed next gate $j'$ and verifier record $\mathbf y$, Bayes' rule gives
 \begin{equation}
\begin{aligned}
 b^-_{iz',j'}&=\sum_z b_{iz}P_i^a(z',j'\mid z,j),\\
 b^+_{iz'}(j',\mathbf y)&=
 \frac{b^-_{iz',j'}P_i^V(\mathbf y\mid z',j',g)}
 {\sum_{\ell} b^-_{i\ell,j'}P_i^V(\mathbf y\mid \ell,j',g)} .
\end{aligned}
\label{eq:belieftransition}
\end{equation} 
The first sum ranges over error locations at the current gate, the second over locations at $j'$. The observed pair has positive predictive probability. At E4 the artifact is unchanged during diagnosis. Its posterior uses the report likelihood $P_i^H(\mathbf d\mid z,\mathbf y)$ conditional on the verifier evidence already available to the specialist. For the illustrative three-stage workflow, Table~\ref{tab:beliefexample} gives a complete update from the prior through a failed test and a subsequent report. The failed test shifts probability toward later-stage errors but leaves both location uncertainty and a 0.068 probability that the artifact is valid. If the controller requests diagnosis and the report points to stage 2, the probability of a stage-2 error rises from 0.338 to 0.737. This increases the relative appeal of checkpoint 1 compared with checkpoint 2, which would retain the likely error: the report changes which completed work is worth preserving.

\begin{table}[htbp]
\centering
\caption{Three-Stage Belief Updates\label{tab:beliefexample}}
{\scriptsize%
\renewcommand{\arraystretch}{1.12}%
\begin{tabularx}{\textwidth}{@{}cL{0.23\textwidth}CCCCC@{}}
\toprule
$z$ & Meaning of $Z_i=z$ & Prior $b_{iz}$ & $\Prb(F\mid z)$ & After test $b_{iz}^{+}$ & $\Prb(D_2\mid z,F)$ & After report $b_{iz}^{H}$\\
\midrule
0 & Artifact is valid & 0.250 & 0.10 & 0.068 & 0.05 & 0.011\\
1 & Earliest error is in stage 1 & 0.300 & 0.20 & 0.162 & 0.10 & 0.051\\
2 & Earliest error is in stage 2 & 0.250 & 0.50 & 0.338 & 0.70 & 0.737\\
3 & Earliest error is in stage 3 & 0.200 & 0.80 & 0.432 & 0.15 & 0.202\\
\bottomrule
\end{tabularx}%
\par}
\smallskip
{\footnotesize\raggedright Notes. $F$ denotes a failed automated test; $D_2$ denotes a report identifying stage 2 as the most likely error location.\par}

\end{table}

For plan $(k,g)$ and remaining slack $\tau$, let $r_i(k,g\mid z,\tau)$ be the probability of timely verified completion from state $z$ under a specified downstream policy. The plan's policy-dependent timely completion probability is
 \begin{equation}
p_i^R(k,g\mid\mathbf b_i,\tau)
=
\sum_{z\in\mathcal Z_i}
b_{iz}r_i(k,g\mid z,\tau).
\label{eq:recoverability}
\end{equation} 
A first repair can fail while leaving enough time for diagnosis and a second repair. Its immediate chance of success and its eventual completion probability in Equation~\eqref{eq:recoverability} can therefore move in different directions. Section~\ref{sec:capacity} studies the checkpoint choice when that later report competes with other workflows for an expert.

\subsection{Markov state and feasible actions}
\label{sec:modelstate}

The state records the information needed to choose an action and predict the next event. For workflow $i$, $\xi_i$ identifies the operating phase, $n_i^F$ counts failed recoveries, $h_i$ records cumulative active expert time, $w_i$ is the current queue age, and $\tau_i=D_i-t$ is the time left before the deadline. The variables $\delta_i^X$, $\delta_i^V$, and $\delta_i^R$ record observable timing information for execution, verification, and recovery. They are residual durations when completion times are known; otherwise they represent elapsed durations together with the conditional remaining-time distributions. The phase distinguishes routine verification from formal certification. Fixed workflow attributes and transition parameters are denoted by $x_i$ and $\vartheta_i$.

The records $\mathbf y_i$ and $\mathbf d_i$ determine which recovery plans are authorized. The controller records the current execution route or recovery checkpoint and route in $a_i^{\rm run}$ and retains this record until the corresponding verification has been processed. Inactive or obsolete entries are set to $\varnothing$. The workflow state is
 \begin{equation}
\mathbf s_i(t)=
\bigl(\xi_i,j_i,c_i,\mathcal K_i(t),\mathbf b_i,\mathbf y_i,\mathbf d_i,n_i^F,h_i,w_i,\tau_i,
\delta_i^X,\delta_i^V,\delta_i^R,{a_i^{\rm run}},x_i,\vartheta_i\bigr).
\label{eq:workflowstate}
\end{equation} 

The shared pool contains $B$ identical experts. A diagnosis has positive service time $S_i^H$ and, once started, continues without interruption. Let $\mathcal Q(t)$ be the waiting set, $\boldsymbol\delta^H(t)$ the observable timing information for busy experts, $\boldsymbol\iota^H(t)$ their workflow assignments, and $\chi(t)$ the arrival-process state and any observed common operating conditions. The full system state is
\begin{equation}
\mathbf S(t)=
\bigl(t,\{\mathbf s_i(t):i\in\mathcal I(t)\},
\mathcal Q(t),\boldsymbol\delta^H(t),\boldsymbol\iota^H(t),\chi(t)\bigr).
\label{eq:systemstate}
\end{equation}
For renewal arrivals $\chi(t)$ includes time since the last arrival. For Poisson arrivals with a known rate, no arrival age is needed.
\begin{assumption}[Workflow dynamics and observations]\label{ass:markov}
Conditional on the observed common operating conditions, initial faults are independent across workflows. A workflow's transitions and signals depend on its error location, recorded attributes and evidence, and the action taken. Their innovations are independent across workflows conditional on these inputs. The arrival process depends on its recorded state. Given the activity, its recorded inputs and elapsed age, residual activity time carries no further information about the hidden fault.
\end{assumption}
The controller can thus update each workflow's belief from its own observations and the common information. An activity still in progress changes its residual-time distribution through its elapsed age.

\paragraph{Feasible event actions.}

At E2, after observing verifier evidence $\mathbf y_i$, the state-dependent action set is
\begin{equation}
\begin{aligned}
\mathcal A_i^{E2}(\mathbf y_i)
={}&\{\mathrm{Advance},\mathrm{Complete},\mathrm{RequestDiagnosis},\mathrm{Abandon}\}\\
&{}\cup\{\mathrm{SubmitCertification}\}\cup\{\mathrm{DirectRework}(k,g):(k,g)\in\mathcal U_i^A(\mathbf y_i)\}.
\end{aligned}
\label{eq:e2actions}
\end{equation}
Advance and Complete require the corresponding verification and certification passes. SubmitCertification requires an eligible artifact and gate, and may precede expert diagnosis. Direct recovery uses preapproved plans. A diagnostic request joins the queue; E3 selects the next case for service.

At E3, let $\mathcal Q^F(t)$ contain the waiting requests that satisfy the service-start rules, have a positive conditional probability of completing diagnosis by $\min\{D_i,T\}$, and leave a feasible postdiagnosis continuation on at least one such outcome. This definition uses the observable state and activity-time kernels, not continuation values. An expected-duration cutoff is a policy approximation and is not part of physical feasibility. The physical action set is
\begin{equation}
\mathcal A^{E3}(\mathbf S(t))=\mathcal Q^F(t)\cup\{\mathrm{Idle}\}.
\label{eq:e3actions}
\end{equation}
Choosing $i$ starts its diagnosis. Section~\ref{sec:dispatch} specifies DRRC's local screen and priority rule. Assignments are sequential when several experts are idle; selecting Idle ends the current assignment sequence.
At E4, after observing expert report $\mathbf d_i$, the action set is
\begin{equation}
\mathcal A_i^{E4}(\mathbf y_i,\mathbf d_i)=
\{\mathrm{SubmitCertification},\mathrm{Abandon}\}
\cup
\{\mathrm{Rework}(k,g):(k,g)\in\mathcal U_i^H(\mathbf y_i,\mathbf d_i)\}.
\label{eq:e4actions}
\end{equation}
Only a passing certification result updates $c_i$, and completion requires all applicable certification gates to have passed. 

\subsection{Objective and exact SMDP}
\label{sec:modelobjective}

If Rework$(k,g)$ is selected, it uses execution-agent time $T_i^R(k,g)$ and returns to an automated gate. Under policy $\pi$, let $N_i(\pi)$ be the total number of expert-diagnosis visits generated by workflow $i$ over its entire time in the system, including any later visits created after failed recovery. Let $S_{ir}^H$ be the active expert service time of the $r$th visit. We define
\[
H_i(\pi)=\sum_{r=1}^{N_i(\pi)}S_{ir}^H
\]
as \emph{lifecycle expert work}: total active expert time across first and repeat visits, excluding calendar waiting time.

Let $Y_i$ indicate valid completion by the deadline, $E_i$ release of an invalid artifact, and $M_i$ a missed deadline or abandonment, including a workflow unfinished at the horizon. An invalid release can occur when an imperfect routine gate suffices for release. When every release requires formal certification of the final artifact, $E_i=0$ by Assumption~\ref{ass:certification}. Let
$A_i=A_i^{\mathrm{initial}}+A_i^{\mathrm{rework}}$ be total automated execution, $T_i^Q$ total expert-queue time, and $\mathcal C_i^R$ noncompute diagnosis and recovery cost, including fixed diagnostic charges, discarded work, restoration, and route switching. The coefficients $c_H$, $c_A$, and $c_W$ are the respective costs per unit of expert service, automated execution, and queue waiting. For a given pool of $B$ experts, expected net operating value includes the expense of every diagnosis performed. Comparisons across pool sizes account separately for the fixed cost of providing that capacity.
\begin{equation}
\max_{\pi}\ 
\E^\pi\!\left[
\sum_i
\left\{
v_iY_i-\ell_i^E E_i-\ell_i^D M_i
-c_HH_i-c_AA_i-c_WT_i^Q-\mathcal C_i^R
\right\}
\right].
\label{eq:objective}
\end{equation}

On the finite horizon $[0,T]$, expected arrivals are finite, absolute payoffs have an integrable bound, and events are nonexplosive. The action sets are finite and the kernels measurable. Simultaneous events follow a fixed order, with deadlines and the horizon treated as events. Idling advances the clock. Let $\sigma_n$ be the $n$th event time, $\mathbf S_n=\mathbf S(\sigma_n)$ the observed system state and $a_n$ its chosen action. The optimal expected net reward from that event satisfies
\begin{equation}
V_n(\mathbf S_n)=
\max_{a\in\mathcal A(\mathbf S_n)}
\left\{
r(\mathbf S_n,a)+
\E\!\left[V_{n+1}(\mathbf S_{n+1})\mid\mathbf S_n,a\right]
\right\},
\label{eq:smdp}
\end{equation}
Here $r(\mathbf S_n,a)$ records the expected reward and costs assigned to the transition ending at the next event or $T$. Full service liabilities for nonpreemptive diagnoses are booked at service start, including service that continues beyond a workflow deadline or the decision horizon. The terminal value at $T$ charges only outcomes not already recorded; unfinished workflows incur their terminal loss once. Other workflows may join the queue while one agent attempts a repair. Section~\ref{sec:capacity} isolates this dependence by holding the focal workflow's evidence and recovery technology fixed while changing the expert capacity available to competing requests.

\section{Diagnostic Information, Expert Capacity, and Recovery}
\label{sec:capacity}\label{sec:recoverydemand}

How does competition for a specialist change the work a workflow preserves? After failed verification (E2), the controller can request diagnosis or select a checkpoint for workflow $A$'s first repair. We compare these choices \emph{before} observing the repair outcome, holding the fault distribution, recovery technology, and incoming workflows fixed as expert capacity changes. This separates the operational response to capacity from improvements in execution technology. The simplified timing below isolates how competition for later advice changes a current recovery decision. We consider a known fault whose location remains uncertain.

\subsection{The value of access to diagnosis after a trial}

Workflow $A$ reaches screening at time 0; a known group of $n$ other workflows arrives at time 2. All have deadline 5 and share $B\geq1$ experts. Diagnosis takes one unit of expert time; repair takes two units of agent time. Workflows can repair in parallel but each undertakes one operation at a time. An immediate repair of $A$ is verified at time 2: success completes it, while failure leaves time for another repair, with or without a report, or abandonment. These durations leave time for diagnosis followed by repair, either immediately or after a failed trial. Immediate diagnosis finishes before the other workflows arrive, whereas diagnosis after failure competes with them. The timing thus isolates the value of preserving access to later advice. Completions precede new starts. Timely completion earns $v_i>0$; waiting, abandonment and expiry carry zero cost.

The route is fixed, and the same checkpoints $\mathcal K_i$ are available before and after diagnosis. Checkpoint $k$ retains stages 1 through $k$ and re-executes the remainder at cost $C_i^R(k)\geq0$; checkpoint 0 rebuilds everything. Fault locations $Z_i\in\{1,\ldots,J_i\}$ are independent with priors $\mathbf b_i$. Repair succeeds precisely when $Z_i>k$; failure preserves the fault and reveals $Z_i\leq k$. Retaining work can therefore save expenditure without saving time. A report identifies the fault at cost $C_i^H=c_i^0+c_H S_i^H>0$, including the visit charge and active service. Suppressing $i$, the last opportunity for recovery yields
 \begin{equation}
\begin{aligned}
Q_{\rm last}^0(\mathbf b)
 &=\max\!\left\{0,\max_{j\in\mathcal K}\left(v\sum_{z>j} b_z-C^R(j)\right)\right\},\qquad
 \Delta_{\rm last}=Q_{\rm last}^{H}-Q_{\rm last}^{0},\\
Q_{\rm last}^H(\mathbf b)
 &=\sum_z b_z\max\!\left\{0,\max_{j\in\mathcal K:\,j<z}(v-C^R(j))\right\}-C^H .
\end{aligned}
\label{eq:capterminal}
\end{equation} 
The gain $\Delta_{\rm last}(\mathbf b)$ measures how much a report improves the choice of a final repair after paying for diagnosis. For the $n$ arrivals at time 2, order $[\Delta_{{\rm last},i}(\mathbf b_i)]_+$ from largest to smallest as $d_{(1)},d_{(2)},\ldots$, padding with zeros, where $[x]_+=\max\{x,0\}$. On their own these workflows earn $W_B=\sum_i Q_{{\rm last},i}^0(\mathbf b_i)+\sum_{j=1}^B d_{(j)}$. If $A$ takes an expert at time 2, they forgo $\theta_B=d_{(B)}$. This is zero when a specialist would otherwise be free or the displaced report has no value.

For $A$, a trial from $k$ fails with probability $q_k=\Prb(Z\leq k)$. On that branch the posterior is $b_z^k=b_z\mathbf 1\{z\leq k\}/q_k$. If $A$ proceeds without a subsequent report, including the option to repair again on this evidence, the trial contributes $U_k=-C^R(k)+(1-q_k)v+q_kQ_{\rm last}^0(\mathbf b^k)$. Terms weighted by $q_k$ are zero when failure is impossible.

\begin{lemma}[Value of a recovery attempt]\label{lem:capvalue}
A first recovery from $k$ contributes
 \begin{equation}
Q_k(B)=U_k+q_k[\Delta_{\rm last}(\mathbf b^k)-\theta_B]_+
\label{eq:capoption}
\end{equation} 
to system value over and above $W_B$. Expected net operating value, before the fixed cost of providing $B$ experts, is
$V(B)=W_B+\max\{0,Q_{\rm last}^H(\mathbf b),\max_k Q_k(B)\}$.
\end{lemma}
The second term in~\eqref{eq:capoption} is the option value of advice after an unsuccessful trial, net of the diagnostic opportunity cost imposed on competing workflows. At equality either allocation is optimal. The checkpoint determines both the probability of returning to screening and the work a report can help preserve. A direct repair may avoid expert service altogether; if it fails, it can change the evidence and the value of requesting advice. Thus even an action that uses no expert now must be evaluated together with the diagnostic demand it may generate later. Appendix~\ref{app:capacityproof} proves the results and treats competing workflows whose characteristics become known only on arrival.

\subsection{Expert capacity and checkpoint choice}

Suppose a trial from $k$ has failed, so the fault lies in the retained prefix. Without a report, a repair from $j<k$ saves $C^R(0)-C^R(j)$ relative to rebuilding, but forfeits expected completion value $v q_j/q_k$. Restarting at or after $k$ would leave the fault in place. The following condition makes a full rebuild from checkpoint 0 optimal on every reachable failure branch.
\begin{condition}[Unaided rebuilding]\label{cond:caprebuild}
Let $\bar k=\max\mathcal K$. The focal workflow has
 \begin{equation}
 v>C^R(0),\qquad vq_j\geq q_{\bar k}\{C^R(0)-C^R(j)\}
 \quad(j\in\mathcal K,\ j<\bar k).
 \label{eq:caprebuild}
\end{equation} 
\end{condition}
The first inequality makes rebuilding preferable to stopping. The second checks the largest possible failure probability $q_{\bar k}$. For each $j<k$ with $q_k>0$, it ensures that the saving from a partial repair cannot cover the completion value at risk. Let $k^*(B)$ denote the largest maximizer of $Q_k(B)$ among the feasible direct plans.
\begin{theorem}[Capacity and the first checkpoint]\label{thm:capcheckpoint}
Suppose Condition~\ref{cond:caprebuild} holds. For $B_2>B_1$, the largest optimal checkpoint and its first-attempt success probability satisfy
 \[
 k^*(B_2)\geq k^*(B_1),\qquad
 1-q_{k^*(B_2)}\leq 1-q_{k^*(B_1)} .
\] 
\end{theorem}
With more experts, the optimal direct repair can retain a longer prefix even though this can lower its chance of immediate success. Failure after retaining several stages can leave the fault location unresolved; a subsequent report then helps remove the fault while preserving earlier work. Under Condition~\ref{cond:caprebuild}, a full rebuild from checkpoint 0 would otherwise be optimal after failure. Greater capacity makes the report less costly to competing workflows, so trying to preserve more work becomes more attractive. When a trial is optimal at both pool sizes, the theorem orders the selected checkpoints under the stated tie rule.

\subsection{Diagnosis timing and total expert work}

Diagnosis before $A$'s first repair displaces no other request; diagnosis after failure may do so. The preferred order therefore depends on the competing report value $\theta_B$.

\begin{proposition}[When to diagnose before repair]\label{prop:timing}
The advantage of trying a repair first, $\max_k Q_k(B)-Q_{\rm last}^H(\mathbf b)$, is nondecreasing in the number of experts $B$. Suppose $Q_{\rm last}^H(\mathbf b)>\max\{0,\max_k U_k\}$ and some $q_k>0$. Define
 \begin{equation}
 \theta^\star=\max_{k:q_k>0}
 \left\{\Delta_{\rm last}(\mathbf b^k)-
      \frac{Q_{\rm last}^H(\mathbf b)-U_k}{q_k}\right\}.
 \label{eq:timingcutoff}
\end{equation} 
The controller strictly prefers to request a report before repair when $\theta_B>\theta^\star$ and to try a repair first when $\theta_B<\theta^\star$. At equality either order is optimal.
\end{proposition}
For a trial from $k$, the option to consult an expert after failure must cover the initial advantage of diagnosis over unaided recovery. Only a fraction $q_k$ of trials reach that decision, which accounts for the denominator in~\eqref{eq:timingcutoff}. The comparison allows any feasible unaided repair after failure.

\begin{corollary}[Capacity expansion can reduce total expert work]\label{cor:capwork}
Suppose the cutoff in Proposition~\ref{prop:timing} satisfies $0<\theta^\star<d_{(n)}$ for the $n\geq1$ competing workflows. Then $A$ obtains a report before repair at $B=n$ and tries a repair first at $B=n+1$. For any optimal initial checkpoint $k$ at the larger pool,
 \[
 V(n+1)>V(n),\qquad H(n)-H(n+1)=1-q_k>0 .
\] 
Here $H(B)$ is expected total active expert time under the indicated optimal policy.
\end{corollary}
The cutoff compares with the least valuable competing report. At $B=n$, a returning workflow would displace that report; at $B=n+1$, both can be served. All $n$ arrivals receive reports at either capacity, but $A$ uses an expert with certainty in the smaller pool and only on failure in the larger one. Concurrent coverage can thus raise net operating value while reducing actual diagnostic work. Section~\ref{sec:actionableinformation} measures both outcomes and considers competing requests also present at the outset.

\section{Structure-Guided Diagnosis--Resource--Recovery Control}
\label{sec:drrc}\label{sec:policy}
In the simplified competition model of Section~\ref{sec:capacity}, the controller knows which workflows will compete for diagnosis after a failed repair. In the general system, arrivals and ongoing repairs can change the queue. DRRC therefore looks ahead a few recovery attempts within each workflow, retaining the option to request advice at later screening decisions as in Lemma~\ref{lem:capvalue}. The structural analysis motivates retaining contingent diagnosis in this local recovery calculation. Predicted waits, a charge $\mu\geq0$ per unit of expert time, and the dispatch rule approximate competition in the general system. The charge represents the opportunity cost of scarce expert work in the local approximation, while $c_H$ records actual service expenditure. The pool of $B$ experts determines which services can actually start.

\subsection{Computing the DRRC decisions}
\label{sec:localvalues}\label{sec:dispatch}\label{sec:priceimplementation}

\emph{Screening (E2).} DRRC first evaluates each feasible direct plan $u=(k,g)$ from the current belief $\mathbf b$. It works backward over at most $d$ repairs, subject to the remaining attempt limit $m$ and time $\tau$ before the deadline or horizon. Diagnosis consumes time but does not reduce $d$ or $m$. At an intermediate gate it uses~\eqref{eq:belieftransition} to average continuation values over the next observations, including a return to execution after a successful repair. For the terminal-gate implementation used in Section~\ref{sec:experiment}, success earns $v_i$ and abandonment has zero payoff. The value of trying $u$ is then
\[
 \widehat W_i^R(u;\mu)
 = v_i p_u-c_u
   +(1-p_u)\widehat V_i(\mathbf b^F_u,\tau-t_u,d-1,m-1;\mu).
\]
Here $c_u$ and $t_u$ are the cost and duration of the repair, and $p_u=\mathbf b\,\widehat{\mathbf p}_u$ averages its estimated success probabilities across fault locations. The matrix $\widehat K_u^{\rm fail}(z,z')$ gives the joint probability of failure and next location $z'$ from $z$. After failure, $\mathbf b^F_u=\mathbf b\widehat K_u^{\rm fail}/(1-p_u)$ when $p_u<1$. The local value $\widehat V_i$ takes the best feasible action at the new screening decision, including requesting a report. Only plans that fit the remaining time are considered. The failure term is zero if $p_u=1$, $d=1$, or $m=1$. Reaching the lookahead limit ends this calculation, while the physical attempt limit governs the workflow.

Equivalently, following the continuation selected at price $\mu$ gives
\begin{equation}
\begin{aligned}
\widehat W_i^R(u;\mu)&=\widehat B_i(u)-\mu\widehat H_i^+(u),\\
\widehat H_i^+(u)&=(1-p_u)\widehat H_i(\mathbf b_u^F,\tau-t_u,d-1,m-1;\mu).
\end{aligned}
\label{eq:recoveryrelaxedvalue}
\end{equation}
Here $\widehat B_i(u)$ is operating payoff net of recovery, diagnostic, and waiting costs, and $\widehat H_i$ tracks expert time under the continuation selected by $\widehat V_i$. Thus $\widehat H_i^+(u)$ measures expert work induced by repair failure, with the same stopping rules as the value recursion. This later service is already priced in $\widehat V_i$. Writing $Q_i^a(\mathbf s_i;\mu)$ for the value of action $a$ followed by the local continuation, E2 records the best action without an immediate request and its value,
\begin{equation}
 Q_i^0(\mathbf s_i;\mu)=
 \max_{a\in\mathcal A_i^{E2}(\mathbf y_i)\setminus\{\mathrm{RequestDiagnosis}\}} Q_i^a(\mathbf s_i;\mu).
 \label{eq:bestnonexpert}
\end{equation}
An immediate request is valued before its report is known. For each possible report $\mathbf d$, DRRC updates the belief and deducts predicted wait $\widehat w_i$ and service time $\widehat s_i$ to obtain the projected E4 state $\mathbf s_i^H(\mathbf d)$. It then averages the best postreport values and deducts diagnostic costs:
\begin{align}
 W_i^H(\mathbf s_i,\mathbf d;\mu)
  &=\max_{a\in\mathcal A_i^{E4}(\mathbf y_i,\mathbf d)}
        Q_i^a(\mathbf s_i^H(\mathbf d);\mu),\label{eq:postreportvalue}\\
 Q_i^H(\mathbf s_i;\mu)
  &=\widehat{\E}_{\mathsf D_i}\bigl[W_i^H(\mathbf s_i,\mathsf D_i;\mu)\bigr]
        -c_i^0-(c_H+\mu)\widehat s_i-c_W\widehat w_i .
 \label{eq:escalationvalue}
\end{align}
The work calculation follows the same choices. Let $\widehat H_i^0$ and $\widehat H_i^{E4}(\mathbf d)$ be future expert time under the actions selected in~\eqref{eq:bestnonexpert} and~\eqref{eq:postreportvalue}, respectively, with their state arguments suppressed. The request entails $\widehat s_i+\widehat{\E}_{\mathsf D_i}[\widehat H_i^{E4}(\mathsf D_i)]$ units of expert time. At a screening decision the recursion closes as
\begin{equation}
(\widehat V_i,\widehat H_i)=
\begin{cases}
\bigl(Q_i^H,\widehat s_i+\widehat{\E}_{\mathsf D_i}[\widehat H_i^{E4}(\mathsf D_i)]\bigr),
 &\text{if a request is feasible and }Q_i^H>Q_i^0,\\
(Q_i^0,\widehat H_i^0),&\text{otherwise}.
\end{cases}
\label{eq:localvaluework}
\end{equation}
Abandonment supplies $(0,0)$. The terminal-gate E4 choice is recovery or stopping; another report requires a failed recovery and a fresh E2 decision. Ties favor the action without a new request.

For a fixed decision time and recovery menu, the report has positive information value precisely when no single recovery plan is optimal for almost every report (Appendix~\ref{app:proofs}). Even an accurate report adds no information value if rebuilding remains optimal wherever the fault lies. The net request advantage $\Delta_i^H=Q_i^H-Q_i^0$ additionally accounts for delay, operating costs, and the charge for expert work.

The wait forecast $\widehat w_i$ scales queued predicted work plus conditional residual service work by a validation-selected coefficient and divides by $B$. The charge is $\mu=\alpha_\mu\max\{0,\rho_0-\rho_\mu\}$, where $\rho_0$ is arrival rate times mean initial service duration divided by $B$, $\rho_\mu\geq0$ is the load baseline, and $\alpha_\mu\geq0$ is selected on validation episodes. Both inputs remain fixed within a local calculation. Policies use mean service durations; the simulation samples actual durations and enforces expiry. Appendices~\ref{app:proofs} and~\ref{app:repro-models} give the auxiliary work-budget interpretation and estimation details.

\emph{Expert dispatch (E3).} Once a specialist is free, waiting requests compete for actual service. DRRC checks which reports can still be used before their deadlines. For each physically feasible request in $\mathcal Q^F(t)$, it recalculates the diagnosis-conditioned value in~\eqref{eq:escalationvalue} with a one-recovery lookahead and zero additional wait. Denote this value by $\widehat D_i^{E3}(t;\mu)$. DRRC's local positive-value screen retains requests in
\begin{equation}
\widehat{\mathcal Q}^{+}(t)=
\{i\in\mathcal Q^{F}(t):\widehat D_i^{E3}(t;\mu)>0\}.
\label{eq:drrceligibility}
\end{equation}
E2 compares requesting with acting without a report. Under the queue commitment, E3 screens the local value of starting service now rather than reopening the request-versus-direct-recovery comparison. Eligible requests are assigned priority levels and scores as follows:
 \begin{align}
L_i&=\max\{0,\min(\bar L,n_i^F)-\lfloor w_i/A_H\rfloor\},\label{eq:level}\\
\Pi_i^{\rm DRRC}&=\frac{(1-\gamma)\widehat p_i^{R,1}v_i+\gamma v_{\rm ref}}
 {\max\{\varepsilon_\Pi,\exp(\delta\widehat s_i)-1\}}\nonumber\\
&\qquad+\eta\left[\frac{w_i-A_H\lfloor w_i/A_H\rfloor}{A_H}
 +\min\!\left\{1,\frac{\widehat s_i+b_H}{\max(\tau_i,\varepsilon_\tau)}\right\}\right].
\label{eq:regimescore}
\end{align} 
Here $w_i$, $\tau_i$, and $n_i^F$ denote queue age, remaining time, and failed recoveries. The proxy $\widehat p_i^{R,1}=\max_u\mathbf b_i\widehat{\mathbf p}_u$ uses all supplied plans before timing restrictions; feasibility is checked separately at dispatch and recovery. The load weight is $\gamma=\operatorname{clip}_{[0,1]}\{(\rho_0-\rho_L)/(\rho_U-\rho_L)\}$, with $0\leq\rho_L<\rho_U$. The nonnegative integer $\bar L$ caps priority levels, $A_H>0$ sets the promotion interval, and $b_H\geq0$ is a deadline buffer. The service penalty rate $\delta>0$ has inverse-time units; $v_{\rm ref}>0$ and $\eta\geq0$ have the payoff units of $v_i$. Positive floors $\varepsilon_\Pi$ and $\varepsilon_\tau$ stabilize the dimensionless and time denominators. The reward--duration ratio is motivated by discounted-reward scheduling in a fixed queue (Appendix~\ref{app:dispatchdetails}, which also gives the fixed settings). Priority levels, aging, and urgency adapt this rule to repeat requests and deadlines. Repeated failures lower priority, while waiting promotes requests toward level zero. Within a level, the score balances recovery reward against service duration, giving shorter service more weight as load rises; its bracketed terms capture aging and deadline urgency. A free expert serves the largest score in the lowest occupied level. Service is nonpreemptive, assignments update the queue sequentially, and experts idle when no eligible request remains.

\emph{Action after diagnosis (E4).} The actual report updates the belief, and the time remaining determines the choice in~\eqref{eq:postreportvalue}: a feasible checkpoint--route plan, certification, or abandonment. Repair returns to verification and, if needed, a fresh E2 comparison. Algorithm~\ref{alg:drrc} summarizes the sequence.

\begin{algorithm}[!htbp]
\caption{DRRC at workflow and expert decision events}
\label{alg:drrc}
{\begin{list}{}{\setlength{\leftmargin}{0pt}\setlength{\itemsep}{2pt}\setlength{\parsep}{0pt}\setlength{\topsep}{0pt}}
\small
\item Fix the estimates of recovery quality and diagnostic duration and the control parameters.
\item \textbf{while} a scheduled event remains through $T$ \textbf{do}
\begin{list}{}{\setlength{\leftmargin}{1.4em}\setlength{\itemsep}{2pt}\setlength{\parsep}{0pt}\setlength{\topsep}{2pt}}
 \item Select the next event and update the clock, activity ages and remaining times. At a common time, process incoming evidence before dispatch.
 \item \textbf{if} a workflow reaches screening (E2) \textbf{then}
\begin{list}{}{\setlength{\leftmargin}{1.4em}\setlength{\itemsep}{2pt}\setlength{\parsep}{0pt}\setlength{\topsep}{2pt}}
   \item Compute its best action without a current request and its value $Q_i^0$. Request diagnosis if feasible and $Q_i^H>Q_i^0$. Otherwise take that best action.
 \end{list}
\item \textbf{else if} a request finds an idle expert or a specialist becomes free (E3) \textbf{then}
\begin{list}{}{\setlength{\leftmargin}{1.4em}\setlength{\itemsep}{2pt}\setlength{\parsep}{0pt}\setlength{\topsep}{2pt}}
   \item Use~\eqref{eq:drrceligibility} to form the eligible queue. Serve the highest-scoring request at the lowest occupied level, update the queue, and repeat while an expert and an eligible request remain.
 \end{list}
\item \textbf{else if} a report arrives (E4) \textbf{then}
\begin{list}{}{\setlength{\leftmargin}{1.4em}\setlength{\itemsep}{2pt}\setlength{\parsep}{0pt}\setlength{\topsep}{2pt}}
   \item Update the belief and use the actual remaining time to choose the action in~\eqref{eq:postreportvalue}.
 \end{list}
\item \textbf{else if} a repair completes (E5) \textbf{then}
\begin{list}{}{\setlength{\leftmargin}{1.4em}\setlength{\itemsep}{2pt}\setlength{\parsep}{0pt}\setlength{\topsep}{2pt}}
   \item Send its output to verification. The verifier result triggers E2.
 \end{list}
\item \textbf{end if}
\end{list}
\item \textbf{end while}
\end{list}}
\end{algorithm}

\subsection{Learning-augmented DRRC}
\label{sec:learning}
DRRC uses population estimates of recovery quality and diagnostic duration. LA-DRRC fits case-specific success probabilities and service times from historical predecision attributes, route outcomes, and diagnostic work, then uses the same E2, E3, and E4 calculations. Quality estimates affect checkpoint--route choices; time estimates affect whether advice leaves time for recovery. The population/contextual distinction concerns the quality and duration inputs used for lookahead control, including the failure kernels for prospective repairs. In the controlled simulations, both policies use the same fitted contextual evidence filter after an observed outcome. The simulations use gradient-boosted predictors and separately validate each controller's price slope, wait multiplier, and lookahead. Section~\ref{sec:componentcontributions} isolates prediction effects at common control settings. Appendix~\ref{app:repro-models} gives fitting details, and Section~\ref{sec:empirical} covers held-out software repositories.

\section{Numerical Evaluation}
\label{sec:experiment}
The controlled experiments examine operating performance under rising demand, changes in recovery decisions after capacity expansion, and the contributions of information and control rules. We first vary arrivals with two experts, then vary capacity for fixed arrivals, and finally compare recovery menus, report quality, controller components, and independent UCT and PPO benchmarks. A software-repair digital twin evaluates the policies using observed recovery outcomes.

\subsection{Experimental design and comparison policies}
\label{sec:closedlooptestbed}
Each episode contains 180 workflows with six stages ($J_i=6$), saved checkpoints $\{0,\ldots,5\}$, and at most four recovery attempts. Each enters screening after its final-stage check detects a fault of unknown location. Two experts supply noisy location reports. The controller selects a targeted or rebuild route from a saved checkpoint; retaining a longer prefix costs less but preserves faults within it. Timely success earns the completion value, while failure updates the belief and returns the workflow to screening if an attempt remains. Started diagnosis occupies its expert for the full service duration, including service after expiry.

An expert report identifies the true error location with probability 0.78. Diagnostic service is lognormal with coefficient of variation (CV) 0.45, and arrivals are Poisson. Workflow values and deadline slacks vary across tasks. Appendix~\ref{app:reproducibility} specifies their distributions together with recovery costs, success probabilities, and failure transitions.

Nominal load measures potential diagnostic workload relative to expert capacity: $\rho_0=\lambda\bar s_{\rm gen}/B$, with arrival rate $\lambda$ and cohort mean initial diagnostic time $\bar s_{\rm gen}$. It counts one potential diagnosis per arrival, before selective requests and repeat visits. With $B=2$, we vary load over 0.50, 0.80, 1.05, and 1.25. At 1.25, diagnosing every arrival once would require 125\% of average capacity. Realized utilization depends on selective requests and repeat visits.

The verifier immediately observes whether a repair succeeds. A report changes the belief about fault location, with the same checkpoint and route menu available on either side of diagnosis. Initial execution costs are sunk. This experiment has no certification charge, false-release loss or abandonment penalty.

Recovery-quality and service-time estimates use 24,000 three-attribute training contexts, 4,000 validation contexts, and 6,000 independent prediction-test contexts. Population means and contextual gradient-boosting models are frozen before policy testing (Appendix~\ref{app:repro-models}).

Table~\ref{tab:policyroles} describes the policies and the role of each comparison. Automated DP optimizes direct-only recovery under its fitted transition model. Section~\ref{sec:runtime} evaluates UCT and PPO alongside both versions of DRRC on common independent episodes.

\begin{table}[htbp]
\centering
\caption{Comparison Policies in the Controlled Experiments\label{tab:policyroles}}
{\small\begin{tabular}{@{}L{0.18\textwidth}L{0.40\textwidth}L{0.34\textwidth}@{}}
\toprule
Policy & Decision rule & Purpose of comparison \\
\midrule
Automated DP & Chooses direct recoveries over all remaining attempts without requesting a report. & Value of expert access beyond optimized automated recovery. \\
Diagnose-all FIFO & Requests a report at every eligible failure, serves requests in arrival order, and uses the report for one-step recovery. & Cost of universal requests under shared capacity. \\
One-step VOI & Requests a report when the value of information for the next recovery covers its cost. & Value of looking beyond the next attempt under selective diagnosis. \\
DRRC & Coordinates requests, dispatch, and recovery using population estimates of recovery quality and service time. & Performance achievable with population inputs. \\
LA-DRRC & Uses the same control structure with contextual quality and service-time estimates. & Performance with contextual inputs. Section~\ref{sec:componentcontributions} isolates their effect at common settings. \\
UCT-256 & Searches over the active workflows using 256 simulations per decision. & Value and decision time of joint online planning. \\
PPO & Learns an actor--critic policy for requests, dispatch, and recovery in simulated episodes. & Value and decision time of learned joint decisions. \\
\bottomrule
\end{tabular}\par}
\smallskip
{\footnotesize\raggedright Notes. Lookahead limits planning depth, not the physical recovery budget.\par}

\end{table}

All policies use the same fitted belief filter and physical service rules. With lognormal service, the remaining time must exceed the shortest feasible recovery duration for diagnosis to have a chance of helping. DRRC and LA-DRRC also apply the value screen in Equation~\eqref{eq:drrceligibility}.

Net value per offered workflow is the primary outcome:
\begin{equation}
\frac{1}{N}\left\{\sum_i v_i\mathbf 1\{i\text{ completes on time}\}
-\sum_{i,r}c_{ir}^{R}-0.02N_H-0.015H-0.015W\right\}.
\label{eq:experimentobjective}
\end{equation}
Costs include recovery expenditure, started diagnoses ($N_H$), active service ($H$), and queue waiting ($W$), including waits that end in expiry. All $N=180$ offers enter the denominator. The charge $\mu$ guides the decisions and is excluded from operating cost.

Validation over 18 combinations selects a two-recovery lookahead and wait multiplier 1.5 for both controllers, with price coefficients 0.09 for DRRC and 0.03 for LA-DRRC. The performance comparisons use 40 new matched episodes per load and pointwise 95\% Student intervals on paired differences. Section~\ref{sec:actionableinformation}'s small instances are solved exactly; the software application separately fits and tunes policies within each training fold.

\subsection{Operating performance under diagnostic load}
\label{sec:operatingperformance}
LA-DRRC earns more per offer than Automated DP at every load (Table~\ref{tab:closedloopmain}). At load 1.25, its net value is 8.7\% higher, while timely completion improves by less than one percentage point; Panel B reports matched uncertainty. The gain therefore reflects the cost of successful recovery as well as completion. Table~\ref{tab:resourceuse} decomposes the costs and expert work.

\begin{table}[htbp]
\centering
\caption{Operating Performance and Paired Comparisons\label{tab:closedloopmain}\label{tab:pairedmain}}
{\small\begin{minipage}{0.92\linewidth}\centering
\emph{Panel A. Mean operating outcomes}\par\smallskip
\begin{tabular}{lrrrrrr}
\toprule
Policy & 0.50 & 0.80 & 1.05 & 1.25 & Complete (\%) & Expert time \\
\midrule
Automated DP & 0.543 & 0.543 & 0.543 & 0.543 & 87.56 & 0.000 \\
Diagnose-all FIFO & 0.554 & 0.444 & 0.290 & 0.181 & 37.19 & 0.897 \\
One-step VOI & 0.597 & 0.586 & 0.576 & 0.566 & 87.49 & 0.749 \\
DRRC & 0.603 & 0.586 & 0.579 & 0.569 & 85.69 & 0.581 \\
LA-DRRC & 0.611 & 0.600 & 0.595 & 0.590 & 88.51 & 0.687 \\
\bottomrule
\end{tabular}
\par\medskip
\emph{Panel B. Paired gains from LA-DRRC at load 1.25}\par\smallskip
\begin{tabular}{lrr}
\toprule
Baseline & Net-value gain [95\% CI] & Completion gain (pp) \\
\midrule
Automated DP & 0.047 [0.039, 0.055] & 0.958 [0.349, 1.567] \\
One-step VOI & 0.024 [0.017, 0.031] & 1.028 [0.313, 1.743] \\
DRRC & 0.021 [0.011, 0.030] & 2.819 [2.080, 3.558] \\
\bottomrule
\end{tabular}
\end{minipage}\par}
\smallskip
{\footnotesize\raggedright Notes. Load columns report net value per offer. Completion and expert service in Panel A are measured at load 1.25; service is per offer. Panel B reports LA-DRRC minus baseline; brackets give 95\% confidence intervals and pp denotes percentage points.\par}

\end{table}

Diagnose-all FIFO completes fewer than two in five workflows at load 1.25: requiring a report at every eligible failure consumes deadline slack even when direct recovery is feasible. Holding universal requests and one-step recovery fixed, shortest-diagnosis dispatch and recovery-$c\mu$ dispatch (highest predicted recovery value per unit of service) raise net value from 0.181 to 0.328 and 0.365, respectively (Appendix Table~\ref{tab:alternativeordering}). Both still send every eligible failure to the queue and trail LA-DRRC's 0.590. At load 1.25, LA-DRRC spends about 19\% less on recovery and makes about 13\% fewer attempts than Automated DP, covering diagnosis and waiting costs. Reports help remove underlying errors while retaining valid earlier stages (Table~\ref{tab:costdecomposition}).

\begin{figure}[htbp]
\centering
\includegraphics[width=\textwidth]{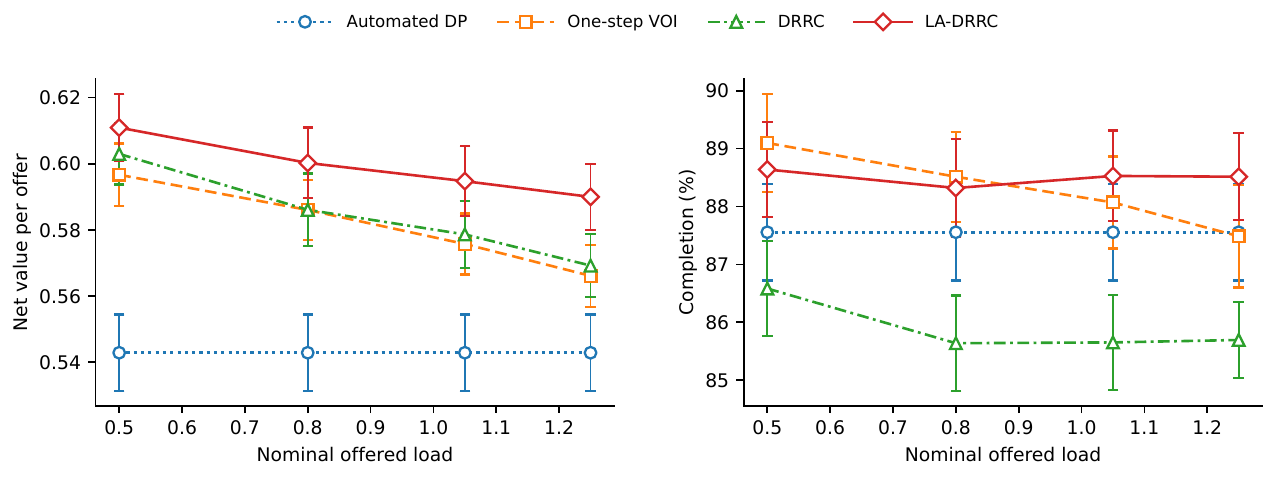}
\caption{Operating performance as nominal diagnostic load increases. Intervals are episode-level 95\% Student intervals.}
\label{fig:closedloopmain}
\end{figure}

\begin{table}[htbp]
\centering
\caption{Operating Costs and Expert Work at Load 1.25\label{tab:resourceuse}\label{tab:costdecomposition}\label{tab:workload}}
{\scriptsize
\begin{minipage}{\textwidth}\centering
\textbf{Panel A. Earned value and costs per offer}\par\smallskip
\begin{tabular}{lrrrrr}
\toprule
Policy & Earned & Recovery & Diagnosis & Waiting & Net \\
\midrule
Automated DP & 0.875 & 0.332 & 0.000 & 0.000 & 0.543 \\
One-step VOI & 0.874 & 0.277 & 0.026 & 0.006 & 0.566 \\
DRRC & 0.857 & 0.267 & 0.020 & 0.001 & 0.569 \\
LA-DRRC & 0.886 & 0.269 & 0.024 & 0.003 & 0.590 \\
\bottomrule
\end{tabular}
\par\medskip
\textbf{Panel B. Diagnostic and recovery work}\par\smallskip
\begin{tabular}{lrrrrrr}
\toprule
Policy & Diagnosed (\%) & Visits & First work & Repeat work & P90 wait & Attempts \\
\midrule
Automated DP & 0.0 & 0.000 & 0.000 & 0.000 & 0.000 & 1.676 \\
One-step VOI & 66.1 & 0.730 & 0.681 & 0.069 & 1.410 & 1.299 \\
DRRC & 51.6 & 0.559 & 0.538 & 0.043 & 0.391 & 1.626 \\
LA-DRRC & 60.9 & 0.677 & 0.622 & 0.065 & 0.802 & 1.453 \\
\bottomrule
\end{tabular}
\end{minipage}\par}
\smallskip
{\footnotesize\raggedright Notes. Diagnosed is the share of offers receiving at least one diagnosis. Visits, attempts, and expert work are per offer; first and repeat work refer to initial and subsequent diagnostic visits. P90 wait averages the episode-level 90th percentiles; Automated DP's P90 wait is set to zero because it has no diagnostic requests. Time is measured in simulation units.\par}

\end{table}

Less expert work is not an objective by itself. Relative to One-step VOI, LA-DRRC uses about 8\% less service and reduces P90 waiting by about 43\% at load 1.25 while improving net value. Relative to separately tuned DRRC, it uses more service and achieves higher completion and value. Anticipating later recovery can avoid low-value requests, while additional intervention is worthwhile when it sufficiently improves recovery. The fixed-parameter comparison below separates contextual prediction from tuning.

\subsection{How the size of the expert pool changes diagnosis and recovery}
\label{sec:actionableinformation}
\emph{Holding the first action fixed.} To measure how much of the return to a larger pool depends on adapting the recovery decision, we keep the arriving workflows fixed. In the instances of Section~\ref{sec:capacity}, we write $V^*(B)$ for the unrestricted optimum $V(B)$ and compute it by exact dynamic programming at $B=1,2$. We then impose $A$'s one-expert optimal first action at two experts and reoptimize all subsequent decisions and all other workflows, obtaining $V^{\rm fixed}(2)$. The same $B$ experts are available throughout each run. In Panel A of Table~\ref{tab:capacitycontrols}, the two workflows arrive at times 0 and 2, earn 10 for completion by time 5, and have identical primitives within each row. A report reveals the fault location in one time unit, and a repair takes two. The timing instance uses fault probabilities $(1/4,1/4,1/2)$, checkpoint costs $(2,1,0)$ in increasing checkpoint order, and diagnosis cost 0.2. For the checkpoint instance the corresponding inputs are $(1/7,1/4,3/28,1/2)$, $(2.6,0.9,0.8,0.1)$, and 0.7.

In the timing instance, the opportunity cost of a time-2 report falls from 1.05 to zero, crossing $\theta^\star=0.20$ in~\eqref{eq:timingcutoff}. With two experts, $A$ tries checkpoint 2 and requests diagnosis only on failure, which occurs with probability one half. Both workflows still finish on time; saving half a 0.2-cost visit raises net value from 18.10 to 18.20. Thus the gain covers at most 0.10 in fixed expansion cost for this cohort. In the checkpoint instance, $A$ repairs first at either capacity; retaining two additional stages saves 0.4893 in expected recovery cost at an additional diagnostic expense of 0.35. Fixing the initial action forfeits the entire operating gain in both instances.

\begin{table}[htbp]
\centering
\caption{Operating Value and Expert Work with One or Two Experts\label{tab:capacitycontrols}}
{\small\setlength{\tabcolsep}{5pt}\begin{tabular}{@{}lllrrrr@{}}
\toprule
& \multicolumn{2}{c}{$A$'s optimal first action} & \multicolumn{3}{c}{Net operating value} & Expert work\\
\cmidrule(lr){2-3}\cmidrule(lr){4-6}
Instance & $B=1$ & $B=2$ & $V^*(1)$ & $V^{\rm fixed}(2)$ & $V^*(2)$ & $H(1)\to H(2)$\\
\midrule
\multicolumn{7}{l}{\textit{Panel A. One workflow arrives initially, one at time 2}}\\
Timing & Diagnose & Repair, $k=2$ & 18.1000 & 18.1000 & 18.2000 & $2\to1.5$\\
Checkpoint & Repair, $k=1$ & Repair, $k=3$ & 17.2964 & 17.2964 & 17.4357 & $1\to1.5$\\
\addlinespace[2pt]
\multicolumn{7}{l}{\textit{Panel B. Three workflows arrive initially, two at time 2}}\\
Timing & Repair, $k=2$ & Diagnose & 44.350 & 44.625 & 44.650 & $3\to5$\\
\bottomrule
\end{tabular}\par}
\smallskip
{\footnotesize\raggedright Notes. Values are exact cohort expectations, net of recovery and diagnostic costs but before fixed staffing costs. $H(B)$ is expected active expert time, equal to expected visits under unit service.\par}

\end{table}

Panel B introduces competition at the initial decision as well. Two workflows arrive with $A$ and must finish by time 4. Their diagnostic costs are 0.1 each. The two workflows arriving at time 2 have deadline 5 and diagnostic costs 0.2 and 1.0. All five share $A$'s value, fault probabilities and repair costs from the timing instance. An immediate report for $A$ now competes with other useful diagnoses. With one expert, $A$ tries a repair first. With two, it seeks advice first and total expert work rises from 3 to 5. Of the 0.300 gain in operating value, 0.275 is available with $A$'s first action fixed. The exact grid confirms that expansion can increase value while either reducing or increasing expert work, and that capturing the full gain can require changing the initial action. Appendix~\ref{app:capacitygrid} reports the complete counts, imperfect-report sensitivity, and treatment of tied optima.

\paragraph{Varying the expert pool at a fixed arrival rate.}
The analytical model isolates competition for diagnosis after a failed trial. In the stochastic system, capacity also relieves competition for immediate diagnosis, as in Panel B of Table~\ref{tab:capacitycontrols}. We match arrivals, attributes, faults, and potential service and recovery draws across pools of one, two, or three experts; the arrival rate gives load 0.80 with two. At the initial screening decision, the automation-first comparator selects its first direct repair using a two-repair lookahead with diagnosis disabled, giving each workflow the same first action across capacities. All subsequent decisions follow LA-DRRC, including the option to request diagnosis. Unlike the exact calculation, which fixes only $A$'s first action, this restriction applies to every arrival.

\begin{table}[htbp]
\centering
\caption{Expert Pool Size and Recovery at a Fixed Arrival Rate\label{tab:fixedarrivals}}
{\small\setlength{\tabcolsep}{4.5pt}\begin{tabular}{crrrrrr}
\toprule
& \multicolumn{2}{c}{Net value per offer} &
\multicolumn{4}{c}{LA-DRRC} \\
\cmidrule(lr){2-3}\cmidrule(l){4-7}
Experts & LA-DRRC & Automation first & Request first (\%) & Expert time & Recovery cost & Complete (\%) \\
\midrule
1 & 0.568 & 0.548 & 44.2 & 0.521 & 0.284 & 87.57 \\
2 & 0.588 & 0.554 & 75.8 & 0.906 & 0.253 & 87.63 \\
3 & 0.599 & 0.555 & 88.7 & 1.092 & 0.240 & 87.82 \\
\bottomrule
\end{tabular}\par}
\smallskip
{\footnotesize\raggedright Notes. Means over 40 matched episodes. Expert time and recovery cost are per offer; request first is the initial-screening request rate. Prediction models and the coefficients relating load to capacity charges and dispatch weights are held fixed across pool sizes.\par}

\end{table}

Expansion from one to three experts raises value by 0.0303 under LA-DRRC and 0.0072 under automation first, a paired difference of 0.0231 (95\% CI $[0.0145,0.0317]$). LA-DRRC's recovery expenditure falls by 0.0447 per offer, covering 0.0194 in additional diagnosis costs; completion changes by only 0.25 percentage points ($[-0.47,0.97]$). Capacity chiefly enables advice before recovery: the initial request rate approximately doubles and expert time rises. All 801 offers starting with direct recovery at both capacities retain the same first checkpoint. The average checkpoint change among direct starters thus reflects selection into direct recovery.

Removing later diagnosis from the initial value calculation, while allowing it after an actual failure, changes first actions on 6.75\% of offers with one expert and 0.94\% with three; the net-value intervals include zero at every pool size (Appendix~\ref{app:fixedarrivals}). With later advice excluded from both initial calculations, permitting an immediate request adds 0.0203 to the one-to-three-expert gain ($[0.0115,0.0292]$). The measurable gain in this population therefore comes chiefly from advice before the initial repair.

\subsection{When diagnostic information improves recovery}
\label{sec:componentcontributions}
With the pool fixed at two experts, we next vary the precision of a report and the recovery choices it can inform. We cross coarse checkpoints $\{0,3\}$ and fine checkpoints $\{0,\ldots,5\}$ with two six-output channels. The location channel identifies the true fault with probability 0.78. The uninformative channel assigns each output probability $1/6$ at every fault location. Diagnosis is available in all four cells. Each cell uses 30 episodes at load 0.80, with common random streams across policies and cells where histories permit; Automated DP and LA-DRRC face the same checkpoint menu.

\begin{table}[htbp]
\centering
\caption{Checkpoint Granularity and Diagnostic Information\label{tab:controlledmechanism}}
{\small\begin{tabular}{lrrrrr}
\toprule
Checkpoints & Report & Auto & LA & Gain [95\% CI] & LA work \\
\midrule
Coarse & Uninformative & 0.504 & 0.503 & -0.001 [-0.005, 0.003] & 0.000 \\
Fine & Uninformative & 0.540 & 0.534 & -0.006 [-0.010, -0.001] & 0.000 \\
Coarse & Location & 0.504 & 0.521 & 0.017 [0.010, 0.024] & 0.533 \\
Fine & Location & 0.540 & 0.598 & 0.058 [0.047, 0.069] & 0.894 \\
\bottomrule
\end{tabular}\par}
\smallskip
{\footnotesize\raggedright Notes. Auto and LA are net value per offer for Automated DP and LA-DRRC. Gain is LA minus Auto, with a paired 95\% interval. LA work is expert service per offer.\par}

\end{table}

With informative reports, LA-DRRC's gain over Automated DP rises from about 3\% with coarse checkpoints to 11\% with fine checkpoints (Table~\ref{tab:controlledmechanism}). Finer checkpoints let the controller act on precise location information while preserving more valid work. Under the fine menu, LA-DRRC discards about 25\% fewer stages per offer, including repeated execution across attempts. With uninformative reports, LA-DRRC requests no diagnosis and slightly trails Automated DP under fine checkpoints because of its shorter recovery horizon. Let $D_{ch,r}$ denote LA-DRRC minus Automated DP in episode $r$, with $c\in\{F,C\}$ for fine or coarse checkpoints and $h\in\{L,U\}$ for location or uninformative reports. Their interaction is
\begin{equation}
\widehat I=\frac{1}{30}\sum_r
\{D_{FL,r}-D_{CL,r}-D_{FU,r}+D_{CU,r}\}=0.046.
\label{eq:checkpointinteraction}
\end{equation}
The paired 95\% interval is $[0.036,0.057]$: finer checkpoints improve LA-DRRC's relative performance more with informative reports. This is a policy-level interaction; checkpoints can also benefit automated recovery, and the comparison charges only for using checkpoints already saved by the platform.

\begin{figure}[htbp]
\centering
\includegraphics[width=\textwidth]{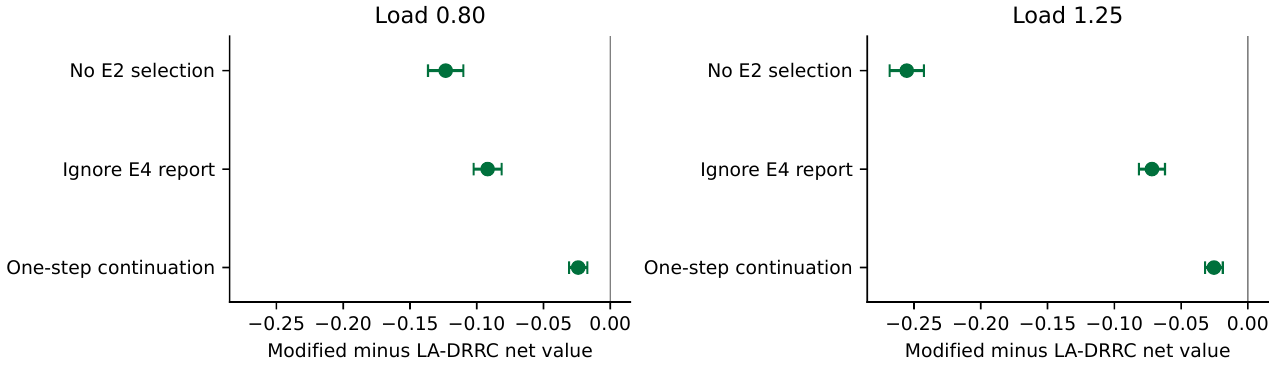}
\caption{Contributions of request selection, report use, and local continuation. Bars show paired 95\% Student confidence intervals.}
\label{fig:ablation}
\end{figure}

We next fix the fine menu, informative report, fitted inputs, and control parameters, changing one LA-DRRC component at a time over 40 new paired episodes at loads 0.80 and 1.25. At load 1.25, universal requests reduce net value by about 44\%, and withholding the report from the immediate recovery choice reduces it by about 12\% (Figure~\ref{fig:ablation}). The former consumes service and deadline slack on avoidable requests; the latter pays for information without using it in the next checkpoint--route choice, although later decisions retain the updated belief. Selective acquisition and effective use of information account for the main gains.

Reducing the lookahead from two recoveries to one lowers net value at load 1.25 by 0.025 per offer (paired 95\% interval $[-0.032,-0.019]$), about 4\%, with no statistically resolved change in completion. Removing only the price on future expert work changes value by $+0.002$ ($[-0.001,0.005]$). Looking ahead improves the choice of recovery plan, while the additional contribution of pricing future service is unresolved. Keeping the E3 eligibility screen and serving eligible requests FIFO changes net value by $+0.003$ ($[-0.001,0.006]$) and raises timely completion by 0.40 percentage points ($[0.05,0.76]$). Thus the detailed ranking has no detectable net-value advantage at this load (Appendix Table~\ref{tab:ablationpaired}). Each intervention changes later trajectories, so its measured effect is conditional on the other rules being in place.

Table~\ref{tab:predictioneffect} in the Appendix holds LA-DRRC's control settings and belief filter fixed. Contextual quality and time inputs together raise net value over population inputs by about 5\% at load 1.25, with quality contributing more individually. Quality estimates guide recovery choices; time estimates indicate whether advice can arrive soon enough. At load 0.80, the time-only interval includes zero, but the combined gain remains positive. Using common settings separates prediction effects from tuning differences in the main DRRC comparison.

Miscalibration can reverse the value of intervention. With actual report accuracy 0.55 but assumed accuracy 0.78, LA-DRRC falls below Automated DP by 0.027 per offer (95\% CI $[-0.037,-0.017]$): overstated accuracy can induce costly requests and overconfident recovery choices. Appendix~\ref{app:supp-sensitivity} examines accuracy, duration, and cost sensitivity.

\subsection{Operating value and decision time}
\label{sec:runtime}
UCT evaluates joint decisions through simulated continuations; PPO learns an actor--critic policy from repeated episodes \citep{schulman2017ppo}. Both and the reference policies face the same 40 independent 180-workflow episodes per load. UCT uses its validation-selected budget of 256 simulations per decision; PPO uses five frozen training seeds, each trained for 524,288 transitions with a validation-selected checkpoint. Table~\ref{tab:benchmarksummary} reports net value from the complete test episodes and decision time measured on common sampled states using the same processor. Appendix~\ref{app:runtime} gives the protocols and load-specific results.

\begin{table}[htbp]
\centering
\caption{Matched Comparison of Operating Value and Decision Time\label{tab:benchmarksummary}}
{\small\begin{tabular}{lrrr}
\toprule
Policy & Net value & Mean time (ms) & Time / LA-DRRC \\
\midrule
Automated DP & 0.541 & 0.729 & 0.99 \\
One-step VOI & 0.577 & 0.614 & 0.84 \\
UCT-256 & 0.452 & 692.939 & 943.96 \\
PPO & 0.521 & 1.143 & 1.56 \\
DRRC & 0.581 & 0.756 & 1.03 \\
LA-DRRC & 0.602 & 0.734 & 1.00 \\
\bottomrule
\end{tabular}\par}
\smallskip
{\footnotesize\raggedright Notes. Net value per offer averages equally across loads; PPO also averages five training seeds. Timing uses 120 common states with equal weights across applicable phases, excluding offline training and shared predictor inference.\par}

\end{table}

\paragraph{Online planning.}
LA-DRRC earns 0.149 more per offer than UCT (95\% CI $[0.142,0.157]$), about 33\% of UCT's value; PPO also outperforms UCT. UCT completes fewer workflows and spends more on recovery than LA-DRRC. Its finite simulation budget spans recovery plans, reports, and competing workflows, while the depth limit and omitted future arrivals restrict its queue forecasts. These limitations may explain its weaker decisions here; DRRC integrates report and recovery outcomes directly within its local recursion.

\paragraph{Neural policy learning.}
LA-DRRC earns 0.080 more per offer than PPO (95\% CI $[0.071,0.089]$ across episodes and training seeds), about 15\% of PPO's value; DRRC also exceeds PPO at all loads. None of the five trained PPO policies requests diagnosis in testing, and diagnostic use declines during training. Immediate request costs precede benefits that require a subsequent informed repair, which may hinder learning as reports become rare. Using only direct recovery, PPO also trails Automated DP by about 3.6\%.

\paragraph{Computational cost.}
UCT takes about 940 times as long per decision as LA-DRRC; PPO takes about 56\% longer. UCT repeatedly simulates the joint system, whereas the local rules reuse a compact recovery calculation. PPO avoids search but incurs feature-construction and neural-inference costs; its training and learning-rate selection require about three recorded CPU hours. Equal phase weights prevent realized diagnostic frequencies from driving the timing comparison.

\subsection{Application to software repair}
\label{sec:empirical}

\paragraph{Repair data and workflow construction.}
The digital twin combines observed direct and reference-assisted repair outcomes from FeedbackEval's SWE-bench records \citep{dai2026feedbackeval} with a simulated shared diagnostic service. Matching the same faulty function across five repair models yields 177 defects from 11 repositories, with published outcomes for every model under direct and report-assisted recovery. Each workflow has one checkpoint---the original faulty function---and may attempt two distinct recovery actions, restoring that input before each. Direct recovery uses test feedback; after diagnosis, the same models can also use the report. The twin retains E2 request selection, two shared nonpreemptive diagnostic servers at E3, and report-dependent E4 recovery. All diagnostic policies use common model-generated reports constructed with reference repairs. Reports become visible only after service; controllers never observe reference repairs or unchosen-action outcomes. Each defect permits at most one diagnostic service, which may follow a failed direct attempt. This single-checkpoint, one-report design evaluates selective access to diagnostic assistance and report-dependent recovery under shared capacity. Checkpoint choice and repeated diagnostic visits are examined in the controlled experiments.

All five direct routes fail on 43 defects, of which 15 admit a successful report-assisted route. Among 885 defect--model pairs, assisted repair succeeds where direct repair fails in 153 pairs; the reverse occurs in 56. These hindsight comparisons characterize available recovery options, whereas the controller must request information before knowing which action will succeed.

\paragraph{Digital-twin operating assumptions.}
The twin adds a common queue, Poisson arrivals, and deadline slack varying with original function size. Timely completion earns one unit; recovery and diagnostic work scale with output and report length, respectively, and incur operating costs. Controllers use work estimates fitted on training data. The twin also allows a queued workflow to return to direct recovery or abandonment if its request fails the dispatch-time screen; diagnosis is disabled at that immediate fallback decision. Appendix~\ref{app:empiricalppo} specifies this rule, timing, and costs.

\begin{figure}[htbp]
\centering
\includegraphics[width=1\textwidth]{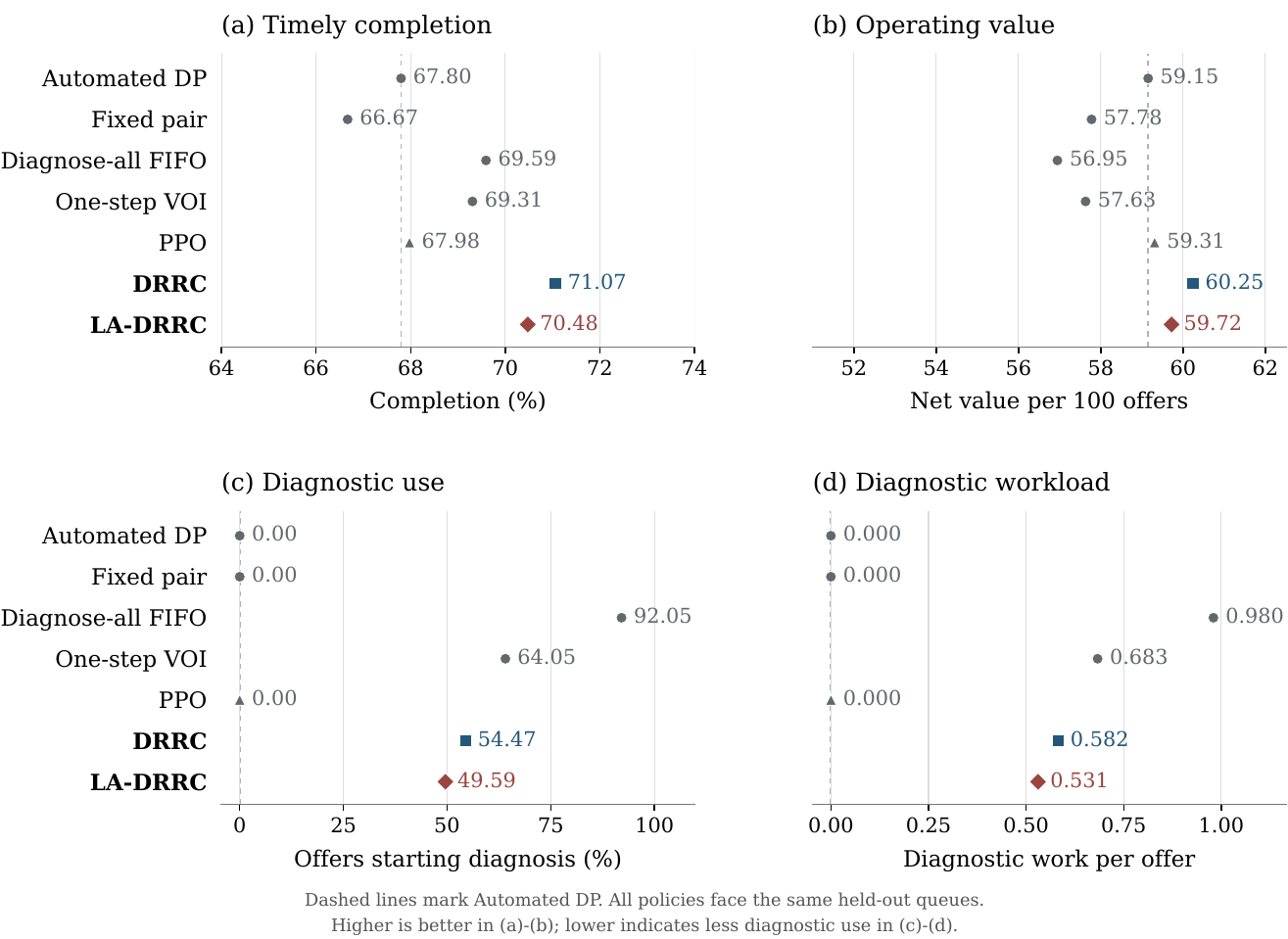}
\caption{Policy performance in the code-repair application. Means weight held-out defects and the four loads equally; PPO also averages five training seeds. Diagnostic work is dimensionless.}
\label{fig:empirical}
\end{figure}

\paragraph{Policy comparison and evaluation design.}
We compare Automated DP, a fixed direct-route pair, Diagnose-all FIFO, One-step VOI, PPO, DRRC, and LA-DRRC. Automated DP optimizes remaining direct attempts and is independently tuned on direct outcomes. The fixed pair, selected by training net value, tries its second route after failure if time permits. In this single-checkpoint application, uncertainty is represented by empirical recovery-response profiles rather than stage-wise fault locations. Each training defect supplies the success outcomes of five direct and five report-assisted routes, together with a report category learned from report text. DRRC starts with equal weights on these profiles; LA-DRRC mixes this distribution with task-specific weights learned from the problem description, faulty code, and test feedback. A failed repair downweights profiles in which that route succeeds, while a returned report reweights profiles by a smoothed report-category likelihood. Posterior-weighted outcomes yield success probabilities for the remaining routes. Report shrinkage moves this likelihood toward an uninformative distribution; the contextual weight controls the mixture of task-specific and population information.

Five repository-disjoint outer folds assess transfer; three inner folds select report shrinkage, contextual weight, and control parameters using training data only. Diagnose-all FIFO and One-step VOI share LA-DRRC's fitted inputs and settings. PPO uses the same fitted response and work estimates but learns joint E2--E4 decisions without DRRC's continuation values or ranking. Inner-fold queue evaluations select its learning rate and budget, followed by five independent training seeds per outer fold (Appendix~\ref{app:empiricalppo}).

\paragraph{Results.}
For each outer fold and each nominal load $0.50$, $0.80$, $1.05$, and $1.25$, the seven policies face 40 matched queues of 120 offers (Figure~\ref{fig:empirical}). Repeated queues reduce simulation noise; the empirical sample remains 177 defects in 11 repositories. Relative to Automated DP, DRRC improves average timely completion by 3.27 percentage points and net value by 1.85\%. Compared with diagnose-all FIFO, it achieves 5.79\% higher net value using 40.57\% less diagnostic work, illustrating the value of selective requests. Contextual augmentation does not improve mean performance over DRRC in this application, although LA-DRRC improves on Automated DP; PPO makes no diagnostic requests. DRRC's mean net-value gain over Automated DP is $0.0109$ per offer, with a repository-clustered 95\% interval of $[-0.0571,0.0354]$. Its corresponding gain over PPO is $0.0094$ ($[-0.0786,0.0395]$, accounting also for training seeds). These intervals leave the value ranking uncertain across projects. The reported gains concern selective access to reference-assisted repair under the specified queue and deadlines.

\FloatBarrier
\section{Managerial Implications and Conclusion}
\label{sec:conclusion}

As automated execution scales, reliable completion depends on coordinating information acquisition with the work that can be preserved and resumed. Expert demand is endogenous: available recovery actions determine the value of advice, while access to experts changes recovery choices and subsequent requests. This feedback links capacity planning to recovery design in AI workflows and other operations combining autonomous execution with shared diagnostic expertise.

Capacity plans should account for when experts must be available and how recovery decisions respond to that availability. In the analytical model, expansion can make it optimal to try recovery first and request advice only on failure, reducing expected visits. In the controlled simulations, it instead attracts more requests before repair and increases expert work. Both responses reflect endogenous diagnostic demand, so historical requests or service hours alone can misstate capacity needs. Expansion can change both diagnosis timing and the work retained; its gain in net operating value must cover the fixed cost of additional capacity. Recovery architecture also matters: checkpoints and alternative routes determine how effectively reports can preserve useful work. Expert attention is most valuable where information changes a feasible downstream action. Evaluating completion together with discarded work, service, and waiting captures these trade-offs more fully than first-attempt success alone. The software application supports selective reference-assisted recovery, although its value advantage remains uncertain across repositories. 

Field evaluation should link chosen restart points to available evidence, expert capacity, later requests, and final outcomes. Natural extensions include choosing checkpoint locations and maintenance expenditure, allocating experts with heterogeneous skills, and adapting decisions as diagnostic quality and recovery performance change. Jointly choosing checkpoint locations and expert capacity would reveal when investment in recoverability substitutes for, or complements, additional staffing. Learning from selectively observed reports and recovery outcomes would further connect day-to-day control to improvements in diagnostic and execution quality.

\clearpage
\appendix
\section*{Appendices}
\addcontentsline{toc}{section}{Appendices}

\makeatletter

\setcounter{section}{0}
\setcounter{subsection}{0}
\setcounter{subsubsection}{0}

\renewcommand{\thesection}{\Alph{section}}
\renewcommand{\thesubsection}{\thesection.\arabic{subsection}}
\renewcommand{\thesubsubsection}{\thesubsection.\arabic{subsubsection}}

\renewcommand{\theHsection}{app.\Alph{section}}
\renewcommand{\theHsubsection}{app.\Alph{section}.\arabic{subsection}}
\renewcommand{\theHsubsubsection}{app.\Alph{section}.\arabic{subsection}.\arabic{subsubsection}}

\setcounter{equation}{0}
\renewcommand{\theequation}{\thesection.\arabic{equation}}
\renewcommand{\theHequation}{app.\Alph{section}.\arabic{equation}}
\@addtoreset{equation}{section}

\setcounter{figure}{0}
\setcounter{table}{0}

\setcounter{algorithm}{0}

\@addtoreset{figure}{section}
\@addtoreset{table}{section}
\@addtoreset{algorithm}{section}

\renewcommand{\thefigure}{\thesection.\arabic{figure}}
\renewcommand{\thetable}{\thesection.\arabic{table}}
\renewcommand{\thealgorithm}{\thesection.\arabic{algorithm}}

\renewcommand{\theHfigure}{app.\Alph{section}.\arabic{figure}}
\renewcommand{\theHtable}{app.\Alph{section}.\arabic{table}}
\renewcommand{\theHalgorithm}{app.\Alph{section}.\arabic{algorithm}}

\makeatother

\section{Proofs and Structural Comparisons}
\subsection{Capacity, checkpoint choice, and service order}\label{app:capacityproof}

\paragraph{Opportunity cost of a diagnosis at time 2 (Lemma~\ref{lem:capvalue}).}
A late arrival has time for a report and a recovery only if its diagnosis starts at time 2. The best allocation adds the largest $B$ positive report gains to the sum of the late workflows' no-report values, with zeros padding the gain list. Adding a candidate with gain $x$ increases this allocation surplus by $[x-d_{(B)}]_+$. Thus, after $A$'s first recovery from $k$ fails, its contribution is $Q_{\rm last}^0(\mathbf b^k)+[\Delta_{\rm last}(\mathbf b^k)-\theta_B]_+$. Weighting by $q_k$, adding the payoff on success and deducting the first recovery cost gives~\eqref{eq:capoption}. A diagnosis before recovery can use an otherwise idle expert and leaves all $B$ experts available when the other workflows arrive.

These alternatives exhaust the initial decision. Before time 2, the first operation can be advanced to time 0 and its outcome held until the original decision time. If $A$ waits until time 2, its incremental value is at most
$\max\{Q_{\rm last}^0(\mathbf b),Q_{\rm last}^H(\mathbf b)-\theta_B\}$.
Its first alternative is attainable by immediate abandonment or direct recovery followed by abandonment on failure, and its second is bounded by the early diagnosis value. A second perfect report supplies no information, and a recovery known to leave the fault in place produces no new evidence. Both waiting and continuous start times are therefore covered.

\paragraph{Ordering the checkpoints (Theorem~\ref{thm:capcheckpoint}).}
For any fault $z$, let $c(z)=\min_{j\in\mathcal K:j<z}C^R(j)$ be the cheapest valid recovery after an exact report. Valid checkpoints accumulate as $z$ increases, so $c(z)$ is nonincreasing. Since $c(z)\le C^R(0)<v$, recovery is profitable after each report. For $j<k$ and $q_k>0$, the success probability from $j$ after failure at $k$ is $1-q_j/q_k$. Equation~\eqref{eq:caprebuild} makes the saved cost no greater than the lost completion value for each such $j$; checkpoints $j\geq k$ retain the fault with certainty. Thus full rebuilding attains $Q_{\rm last}^0(\mathbf b^k)$ after every reachable first failure. Hence
 \[
\Delta_{\rm last}(\mathbf b^k)
  =C^R(0)-C^H-\mathbb E[c(Z)\mid Z\le k],\qquad q_k>0.
\] 
For $j<k$ with $q_j>0$, conditioning on $Z\le k$ adds only locations whose costs $c(z)$ are weakly below every cost for $z\le j$. The conditional mean weakly decreases and the report gain weakly increases. For $B_2>B_1$, the incremental value for checkpoint $k$ is consequently
 \[
 Q_k(B_2)-Q_k(B_1)
 =q_k\int_{\theta_{B_2}}^{\theta_{B_1}}
     \mathbf 1\{x<\Delta_{\rm last}(\mathbf b^k)\}\,dx.
\] 
The integral and $q_k$ are both nonnegative and nondecreasing in $k$; a zero-failure checkpoint gains zero. The maximizer comparison follows from increasing differences \citep{milgrom1994monotone}. Let $k_r$ be the largest maximizer at $B_r$. If $k_2<k_1$, increasing differences gives
$Q_{k_1}(B_2)-Q_{k_2}(B_2)\ge Q_{k_1}(B_1)-Q_{k_2}(B_1)\ge0$.
Then $k_1$ also maximizes value at $B_2$, contradicting the choice of $k_2$. First-attempt success is $1-q_k$, which yields the second comparison.

The same integral orders the checkpoints whenever $[\Delta_{\rm last}(\mathbf b^k)]_+$ is nondecreasing over $k$ with $q_k>0$. Condition~\ref{cond:caprebuild} guarantees this property. To see the role of the unaided fallback, take $v=10$, $\mathbf b=(2,1,4,2)/9$, checkpoint costs $(9,2,1,0)$, diagnosis cost $3/2$, and one competing workflow with a net report gain of 1. After failure at checkpoint 3, restarting at 1 is preferable to rebuilding. The report gains after trials at 2 and 3 are $1/2$ and $-1/14$, respectively. As $\theta_B$ falls from 1 to 0, Equation~\eqref{eq:capoption} moves the optimal checkpoint from 3 to 2.

\paragraph{Timing and total expert work.}
For Proposition~\ref{prop:timing}, $\theta_B=d_{(B)}$ is nonincreasing in $B$. By~\eqref{eq:capoption} each $Q_k(B)$ is therefore nondecreasing, whereas $Q_{\rm last}^H(\mathbf b)$ is independent of $B$. For the cutoff, set $g_k=Q_{\rm last}^H(\mathbf b)-U_k>0$. A zero-failure checkpoint cannot match diagnosis. For $q_k>0$, Lemma~\ref{lem:capvalue} gives $Q_k(B)>Q_{\rm last}^H(\mathbf b)$ precisely when $\theta_B<\Delta_{\rm last}(\mathbf b^k)-g_k/q_k$. Taking the maximum over $k$ gives the stated cutoff and indifference at equality.

For Corollary~\ref{cor:capwork}, $\theta_n=d_{(n)}>\theta^\star>\theta_{n+1}=0$. The cutoff therefore gives a strict switch in order and
$Q_k(n+1)>Q_{\rm last}^H(\mathbf b)>Q_k(n)\geq U_k$
for any optimal trial at the larger pool. Thus $q_k>0$ and a subsequent report has strictly positive value at $B=n+1$, where $\theta_{n+1}=0$. If $q_k=1$, failure leaves the prior unchanged, and trying first followed by diagnosis earns $Q_{\rm last}^H(\mathbf b)-C^R(k)\leq Q_{\rm last}^H(\mathbf b)$, a contradiction. All $n$ later jobs are diagnosed under either allocation, and $W_{n+1}=W_n$. Substitution into Lemma~\ref{lem:capvalue} proves the strict increase in $V$. The remaining visit is certain at $B=n$ and occurs with probability $q_k$ at $B=n+1$. Since each visit lasts one unit, $H(n)-H(n+1)=1-q_k>0$.

\emph{Unobserved late arrivals.}
In Section~\ref{sec:capacity}, suppose the later arrivals are independent of $A$ and their attributes are revealed only at time 2. Assume integrable system payoffs. A trial now earns $U_k+q_k\mathbb E[(\Delta_{\rm last}(\mathbf b^k)-\theta_B)_+]$. Waiting to see the batch yields at most $\max\{Q_{\rm last}^0(\mathbf b),Q_{\rm last}^H(\mathbf b)\}$ in addition to its realized $W_B$, and this bound is attainable before time 2. Thus the initial action maximizes expected plan values. For each realization, the value displaced by a report is nonincreasing in $B$, preserving the timing comparison; taking expectations of the increasing-differences inequality preserves the checkpoint ordering under~\eqref{eq:caprebuild}.
\par

\subsection{Report value and expert work in local decisions}
\label{app:proofs}

\paragraph{Recovery menus and the value of timely diagnosis.}\label{app:checkpointinformation}
For a fixed menu $\mathcal U_m$, write $M_m(b)=\max_{u\in\mathcal U_m}\sum_z b_zr_u(z)$ and $I_m=\E[M_m(b^{\mathsf D})]-M_m(b)$, where $b^{\mathsf D}$ is the posterior after a report. If $u^*$ maximizes payoff at the prior $b$, then $\E[b^{\mathsf D}]=b$ gives
\[
 I_m=\E\!\left[M_m(b^{\mathsf D})-\sum_z b_z^{\mathsf D}r_{u^*}(z)\right].
\]
A report has positive decision value exactly when no single recovery plan is optimal for almost every report. Expanding the menu can improve both the best unaided plan and the choice made after diagnosis. Section~\ref{sec:componentcontributions} measures their combined effect by varying the menu and the information in the report.

\phantomsection\label{app:boundaryproof}
To see how this gain determines the wait worth incurring for a report, suppose two repairs have payoff differences $d_z=u_1(z)-u_2(z)$ with $d_1>0>d_2$. Write $p=\Prb(Z=1)$. Without a report, the operator prefers repair 1 precisely when $p>-d_2/(d_1-d_2)$. Relative to choosing repair 2 in both states, a perfect report earns $pd_1$ by switching to repair 1 in state 1. With only the prior, the best repair adds $\max\{0,pd_1+(1-p)d_2\}$. The gain from observing the location is therefore $\mathcal V^D(p)=\min\{pd_1,(1-p)(-d_2)\}$.

\begin{figure}[htbp]
\centering
\includegraphics[width=.53\textwidth]{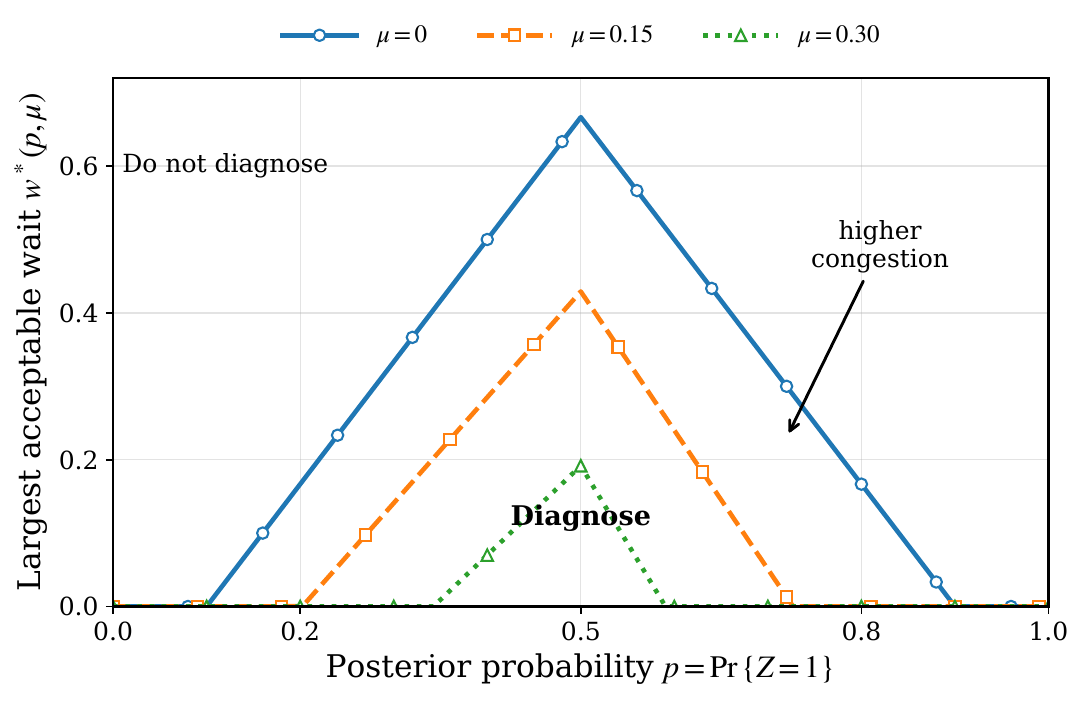}
\caption{Maximum acceptable wait for a diagnostic report in the two-state example.}
\label{fig:reviewboundary}
\end{figure}

Under a fixed diagnostic protocol, let $s_0$ be mean initial service time. In state $z$ each recovery generates a repeat visit with probability $q_z<1$ and mean repeat service $s_z$. Without a deadline or attempt cutoff for these repeat visits, expected later work is $h_z=q_zs_z/(1-q_z)$. The alternative without a report consumes no expert work. Let $K(p)$ be the protocol's expected incremental operating cost at prior $p$, including visit fees, active service on initial and repeat visits, and any value lost during service. Holding the continuation fixed, an initial queue wait of $w$ adds $c_Ww$ to this cost. With $c_W>0$ and a work charge $\mu$, the protocol is worthwhile at a feasible wait $w\geq0$ when
 \begin{equation}
w\leq w^*(p,\mu)=\frac{\mathcal V^D(p)-K(p)-\mu\{s_0+ph_1+(1-p)h_2\}}{c_W}.
\label{eq:waitboundary}
\end{equation} 
Figure~\ref{fig:reviewboundary} illustrates the trade-off with $d_1=-d_2=1$, $K(p)\equiv0.10$, $c_W=0.60$, $s_0=0.20$, $h_1=1.20$, and $h_2=0.30$. At extreme beliefs the preferred repair is already clear. Intermediate beliefs leave scope for an informative report even after some wait.

\paragraph{Pricing expert work over a recovery path.}\label{app:lifecyclecomparison}\label{app:priceinterpretation}
A budget on total expert work gives the charge $\mu$ in Equation~\eqref{eq:recoveryrelaxedvalue} a common unit across workflows. Let $p_{i\pi}$ be a distribution over workflow $i$'s feasible continuation policies, each evaluated under fixed local waits and transitions. Policy $\pi$ earns expected net value $B_{i\pi}$ and uses expected expert work $H_{i\pi}$. Give the planner an aggregate budget $C$ on $\sum_{i,\pi}p_{i\pi}H_{i\pi}$. Relaxing this one budget with multiplier $\mu\geq0$ gives the dual function
\[
 g(\mu)=\mu C+\sum_i\max_{\pi}\{B_{i\pi}-\mu H_{i\pi}\}.
\]
Weak duality bounds the value of the auxiliary expected-work problem. Physical feasibility in the actual queue requires at most $B$ simultaneous nonpreemptive services and depends on their start and finish times. DRRC uses $\mu$ inside the local values while retaining that physical constraint in dispatch.

\section{Model and Policy Implementation}
\label{app:reproducibility}

\subsection{Workflow timing, information, and synthetic inputs}
\label{app:repro-controlled}\label{app:timingbelief}\label{app:requestexpectation}

An observation can change the expected time to completion even when it leaves the error belief unchanged. Under the model's conditional independence of service time and error location, a diagnosis still in progress updates the residual service law using its elapsed age. The expert's report updates the error belief when it arrives. This distinction matters when a queued request might expire before the operator can use the report. Let $\bar\tau_i=\min(D_i,T)-t$, let $\widetilde W_i$ denote the predicted wait until service, and write $I_i^{\rm start}=\mathbf 1\{\widetilde W_i<\bar\tau_i\}$ and $I_i^{\rm report}=\mathbf 1\{\widetilde W_i+S_i^H\le\bar\tau_i\}$. If $\mathcal W_i^H$ is the E4 continuation after the realized report and timing, a request evaluated under the joint timing law has value
\begin{equation}
\E\!\left[-I_i^{\rm start}\{c_i^0+(c_H+\mu)S_i^H\}-c_W\min\{\widetilde W_i,\bar\tau_i\}
+I_i^{\rm report}\mathcal W_i^H-(1-I_i^{\rm report})\ell_i^D\right].
\label{eq:requestexpiry}
\end{equation}
An expired request that never starts incurs no service expense. A diagnosis that starts before expiry occupies an expert for its full duration. The simulator uses realized times; the policies use mean-time forecasts in their continuation calculations.

\paragraph{Synthetic workflow population.}

Contexts are independent uniform draws on $[-1,1]^3$, workflow values on $[0.65,1.35]$, and deadline slacks on $[5,10.5]$. Initial beliefs follow
\[
b\sim{\rm Dirichlet}\{20(0.10,0.14,0.19,0.22,0.20,0.15)\};
\]
the initial error $z\in\{1,\ldots,6\}$ is drawn from $b$. The true route-quality parameters and mean initial diagnostic time are
\begin{align*}
q_T(x)&=\operatorname{clip}_{[0.42,0.94]}\{0.72+0.20\sin(3x_2)+0.11x_1x_3\},\\
q_R(x)&=\operatorname{clip}_{[0.45,0.90]}\{0.68+0.15\cos(2.5x_1)-0.12x_2x_3\},\\
s(x)&=0.65+0.12(x_1+1)+0.42x_2^2+0.25|x_3|.
\end{align*}
These specified functions generate nonlinear heterogeneity in the synthetic workflow population. Gamma arrival gaps have the stated CV and mean $\bar s_{\rm gen}/(\rho_0B)$; CV one gives exponential gaps. Episodes are transient, with no warm-up.

For checkpoint $k$, targeted repair takes $0.5+0.48(6-k)$ time units and rebuild takes $1+0.62(6-k)$. Recovery cost is 0.025 times reference duration plus $0.03(6-k)$ and an additional 0.08 for rebuild. Plans must fit within remaining slack. If $k\ge z$, the retained prefix is invalid: success probability is zero and the error persists. If $k<z$, targeted success is $\max\{0.40,q_T(x)-0.07(z-1-k)\}$ and rebuild success is $q_R(x)$. After failed eligible repair, the reference kernel draws the error uniformly from $\{k+1,\ldots,6\}$. The relative reliability of rebuild and targeted repair depends on context and rollback distance.

The location-report matrix has diagonal $a$ and off-diagonal $(1-a)/5$. The uninformative matrix has all entries $1/6$. Reports update the belief by Bayes' rule; failed plan $u$ updates it by $b'\propto b\widehat K_u^{\rm fail}$ using the common fitted contextual kernel. Controllers do not observe true errors, true quality functions, future arrivals, or future outcomes. Misspecification tests separate the physical report matrix and service multiplier from their assumed values. In the history extension, the known multiplier $1+gf$ enters both physical service and future local diagnosis calculations.

Services are lognormal and nonpreemptive. Queued, diagnosed, or recovering workflows can expire; started expert work is charged in full. Failed fourth attempts exit. All policies use the same physical feasibility constraints. Service and waiting costs enter the objective when incurred.

\subsection{Prediction, continuation, and dispatch}
\label{app:repro-models}\label{app:repro-simulation}

The training, validation, and independent prediction-test seeds are 610001, 610002, and 610003. Both routes are measured at each training context, with 20 Bernoulli trials per quality target and eight service realizations per time target. The population estimator uses the training mean. The contextual estimator uses histogram gradient boosting with 150 iterations, learning rate 0.06, regularization 10, and leaf count selected from $\{7,15,31\}$ by validation squared error. Table~\ref{tab:predictionmetrics} reports prediction errors before the fitted estimates enter the control comparison.

\begin{table}[htbp]
\centering
\caption{Independent Synthetic Prediction Test\label{tab:predictionmetrics}}
{\small\begin{tabular}{llrr}
\toprule
Estimator & Target & Observed-target MSE & Conditional-mean MSE\\
\midrule
Base & Targeted quality & 0.03054 & 0.021717\\
Base & Rebuild quality & 0.01913 & 0.009273\\
Base & Diagnostic time & 0.05239 & 0.025552\\
Learned & Targeted quality & 0.00930 & 0.000167 \\
Learned & Rebuild quality & 0.00993 & 0.000181 \\
Learned & Diagnostic time & 0.02801 & 0.000381 \\
\bottomrule
\end{tabular}\par}
\smallskip
{\footnotesize\raggedright Notes. Observed-target MSE includes trial noise. Conditional-mean MSE uses true generating functions for evaluation only.\par}

\end{table}

The controlled experiments use the following fixed settings for the charge and dispatch rule in Section~\ref{sec:localvalues}:
\begin{align*}
\rho_\mu&=0.50,&\rho_L&=0.65,&\rho_U&=1.30,&\bar L&=3,\\
A_H&=2.5\widehat{\bar s},&b_H&=1.5\widehat{\bar s},&\delta&=0.10,&\eta&=0.10,\\
\varepsilon_\Pi&=10^{-8},&\varepsilon_\tau&=0.01,&v_{\rm ref}&=1.
\end{align*}
Here $\widehat{\bar s}$ is the training-population mean service estimate. Rewards and score weights use the fixed payoff units of the controlled experiments. The endpoints $\rho_L$ and $\rho_U$ specify where the dispatch interpolation reaches zero and one. These implementation settings remain fixed across policy validation and testing. Validation selects the price slope, wait multiplier, and recovery lookahead as described below.

Policy tuning crosses the price coefficients $\{0,0.03,0.09\}$ with wait multipliers $\{0.5,1,1.5\}$ and lookahead depths $\{1,2\}$, using four common validation seeds (62000000--62000003) at every load. The price coefficients are 0.09 for DRRC and 0.03 for LA-DRRC. Both select wait multiplier 1.5 and depth two. The nine simulation test groups use disjoint seeds starting at $73000000+10000s$, with $s=0,\ldots,8$ indexing the main, component, checkpoint--report, accuracy--service, direct-cost, repeat-service, misspecification, one-factor, and joint-change studies.

\paragraph{Implementation details.}
\label{app:localrecursion}\label{app:dispatchdetails}
Section~\ref{sec:localvalues} gives the value--work recursion and dispatch rule. Recovery ties are resolved by fixed plan order. If no recovery fits after a report, its continuation is zero. All report channels are evaluated, including the uninformative channel. In the history extension, future predicted diagnostic duration changes with failure count. Automated DP excludes diagnosis and uses every remaining attempt. The ablation without the future-work price adds $\mu\widehat H_i^+(u)$ back at each recovery comparison.

The service-time term in~\eqref{eq:regimescore} has a scheduling interpretation. In a fixed queue with deterministic service times $s_i$, completion times $C_i$, rewards $R_i$, and discount rate $\delta>0$, pairwise interchange orders jobs by $R_i/(\exp(\delta s_i)-1)$ to maximize $\sum_iR_i e^{-\delta C_i}$. This comparison motivates the reward and service-time terms; DRRC adds priority levels, aging, and urgency to address repeat visits and deadlines.

\section{Experimental Comparisons and Statistical Inference}
\label{app:suppnumerics}

\subsection{Controlled comparisons and operating conditions}
\label{app:supp-prediction}\label{app:supp-ablation}\label{app:informationuncertainty}\label{app:alternativeordering}\label{app:supp-sensitivity}\label{sec:sensitivitymain}\label{app:historystress}\label{app:estimateerror}

\paragraph{Prediction effects at common control settings.}
Table~\ref{tab:predictioneffect} decomposes the effects of contextual quality and time estimates under identical control parameters, as discussed in Section~\ref{sec:componentcontributions}.
\begin{table}[htbp]
\centering
\caption{Local Prediction Effects with Common Parameters\label{tab:predictioneffect}}
{\small\begin{tabular}{lrrrr}
\toprule
Load & Local inputs & Net value & Complete (\%) & Gain over population \\
\midrule
0.80 & Population inputs & 0.584 & 85.35 & 0.000 [0.000, 0.000] \\
0.80 & Quality only & 0.596 & 88.10 & 0.012 [0.005, 0.019] \\
0.80 & Time only & 0.586 & 85.50 & 0.002 [-0.002, 0.006] \\
0.80 & LA-DRRC & 0.599 & 88.56 & 0.016 [0.009, 0.022] \\
1.25 & Population inputs & 0.561 & 84.08 & 0.000 [0.000, 0.000] \\
1.25 & Quality only & 0.581 & 87.83 & 0.020 [0.011, 0.028] \\
1.25 & Time only & 0.569 & 84.83 & 0.008 [0.002, 0.013] \\
1.25 & LA-DRRC & 0.587 & 88.21 & 0.026 [0.019, 0.033] \\
\bottomrule
\end{tabular}\par}
\smallskip
{\footnotesize\raggedright Notes. Population inputs use population estimates for quality and time; Quality only and Time only replace the indicated input with contextual estimates. Each load uses 40 matched episodes. Gains are relative to population inputs, with paired 95\% Student intervals.\par}

\end{table}

For a matched pair of policies, let $D_r$ be their difference on episode $r$. The pointwise 95\% interval is $\bar D\pm t_{0.975,R-1}s_D/\sqrt R$, conditional on the fitted models and control parameters. Main and component comparisons use $R=40$ independent episodes. The checkpoint--report interaction uses 30 episodes. Each episode supplies all four menu--channel outcomes, so the four-cell contrast in Equation~\eqref{eq:checkpointinteraction} is formed within episode before computing its interval. Operating-condition grids use 20 episodes per cell, except the joint stress test, which uses 30 per load.

\paragraph{Grid of first decisions under capacity expansion.}\label{app:capacitygrid}

Across 11,340 exact combinations of fault probabilities, checkpoint savings and diagnostic costs, 4,347 switch strictly from diagnosis first to recovery first. In every one, net operating value rises and all two-expert optimal policies use less expert time than any one-expert optimum. In 1,733 other combinations all optimal policies use more expert time after expansion. In 326 combinations recovery comes first at both capacities and $A$ strictly prefers a later checkpoint with two experts. Of the 6,564 combinations with a positive value gain, 4,759 require a change of initial action to capture it in full. This grid uses exact reports. With imperfect reports, expected visits can rise or fall with expansion. The accuracy comparison below accounts for tied optima.
The focal workflow has four possible error locations. We enumerate all 35 positive priors with probabilities in eighths. Each difference $C^R(k-1)-C^R(k)$ independently takes values in $\{0.5,1.5,3\}$ for $k=1,2,3$, with $C^R(3)=0$. The focal diagnostic cost is 0.25, 0.75, 1.25, or 1.75. The workflow arriving at time 2 has two equally likely error locations and recovery costs $(4,0)$. Its diagnostic cost is 0.25, 1, or 1.75, so assigning the sole expert to $A$ at time 2 forgoes a gain of 1.75, 1, or 0.25. Both workflows have completion value 10. We compare one and two experts in every cell with rational arithmetic, retaining the full set of optimal first actions. When fixing the first action, we give the two-expert system the most valuable member of the one-expert optimal set and reoptimize all continuations. In 4,679 cells this restriction removes the entire positive gain from expansion. Counts describe this grid. The full-rebuild condition holds in 8,604 cells, including 206 strict shifts to a later checkpoint; the remaining 2,736 include 120 such shifts. There is no strict shift to an earlier checkpoint on this grid. The optimal first action is unique and unchanged in 5,784 cells. At least one capacity has tied first actions in the remaining 883 cells outside the strict changes and unchanged cases. For each capacity we compute the minimum and maximum expected visits over all value-optimal policies. A strict increase requires the two-expert minimum to exceed the one-expert maximum. The reverse comparison defines a strict decrease. Evaluating~\eqref{eq:timingcutoff} from the primitives picks out exactly the 4,347 strict changes from diagnosis first to recovery first in this design.

The report-accuracy check uses the two-workflow timing instance in Table~\ref{tab:capacitycontrols}. Give each wrong location probability $(1-a)/2$ when the report is correct with probability $a$. At $a=0.90$, a second expert raises expected visits from 1 to 1.5. At $a=0.95$, visits fall from 2 to 1.5, including when each repair costs an additional 0.05. At $a=0.92$ with the original costs, equally valuable one-expert policies use 1 or 2 expected visits. Randomization spans $[1,2]$, so the direction of change to 1.5 depends on how the original tie is resolved.

\paragraph{Expert pool size at fixed arrivals.}\label{app:fixedarrivals}
The 40 new episode seeds are 91260000--91260039. For each seed the product $B\rho_0$ is held at 1.6, giving identical arrival times across the three pool sizes and nominal loads 1.60, 0.80, and $1.6/3$. Contexts, initial beliefs, faults, deadlines, and values match by workflow. The simulator keys recovery, report, and service draws by episode, workflow, event type, and attempt or visit count. Different actions can thus have different outcomes under common draws. In the automation-first arm the initial plan uses the two-recovery local calculation with diagnosis disabled. Subsequent physical eligibility, dispatch, and recovery follow LA-DRRC. In a second arm the first E2 comparison allows a report now but sets future diagnostic opportunities to zero inside that calculation; all later decisions use ordinary LA-DRRC. Differences in first actions are matched by arriving workflow and include the effects of the resulting queue histories.

The LA-DRRC minus second-arm differences in value per offer are $-0.0022$ $[-0.0051,0.0006]$, $-0.0004$ $[-0.0017,0.0009]$, and $0.0005$ $[-0.0005,0.0015]$ for one, two, and three experts, respectively. Intervals are paired 95\% Student intervals over episodes. Across the same capacities, average first-attempt success under LA-DRRC is 60.95\%, 63.91\%, and 64.70\%, including recoveries chosen after a report. Among direct starters, the average first checkpoint is 1.94, 1.80, and 1.70. For automation first it is 2.07 at every capacity, with a first success rate of 55.04\%. DRRC with population estimates attains values 0.522, 0.572, and 0.589. As a check on the load input used by the controller, fixing both the capacity charge and load-dependent dispatch weights at their two-expert values while updating the wait estimate and physical service capacity gives LA-DRRC values 0.559, 0.588, and 0.599. The one-to-three gain is 0.0403 $[0.0317,0.0489]$.

\paragraph{Acquiring and using diagnostic information.}
Forcing diagnosis whenever it is feasible changes E2 request selection. At E3 the FIFO variant keeps the same positive-value eligibility screen and serves its oldest eligible request. Suppressing the report at E4 uses the saved pre-report belief for the immediate recovery and retains the posterior for subsequent decisions. One-step continuation removes one recovery step from the planning horizon. Removing the charge on future expert work leaves the charge on current service in place. Each intervention uses the original validation settings. Table~\ref{tab:ablationpaired} shows how each affects operating value and timely completion; the separate increment due to future-work pricing is unresolved at both loads.

\begin{table}[htbp]
\centering
\caption{Paired Changes after Modifying LA-DRRC\label{tab:ablationpaired}}
{\scriptsize\begin{tabular}{lrrr}
\toprule
Load & Change & Net-value change [95\% CI] & Completion change (pp) \\
\midrule
0.80 & No E2 selection & -0.123 [-0.136, -0.110] & -18.181 [-19.777, -16.584] \\
0.80 & FIFO among E3-eligible requests & $-0.001$ [$-0.003$, $0.001$] & $0.028$ [$-0.152$, $0.207$] \\
0.80 & Ignore E4 report & -0.092 [-0.102, -0.081] & -4.361 [-5.201, -3.521] \\
0.80 & One-step continuation & -0.024 [-0.031, -0.017] & -0.611 [-1.257, 0.035] \\
0.80 & No future-work price & -0.000 [-0.001, 0.001] & -0.056 [-0.144, 0.033] \\
1.25 & No E2 selection & -0.255 [-0.268, -0.243] & -34.542 [-36.031, -33.052] \\
1.25 & FIFO among E3-eligible requests & $0.003$ [$-0.001$, $0.006$] & $0.403$ [$0.050$, $0.756$] \\
1.25 & Ignore E4 report & -0.072 [-0.082, -0.062] & -3.333 [-4.145, -2.522] \\
1.25 & One-step continuation & -0.025 [-0.032, -0.019] & -0.542 [-1.134, 0.050] \\
1.25 & No future-work price & 0.002 [-0.001, 0.005] & 0.181 [-0.091, 0.452] \\
\bottomrule
\end{tabular}\par}
\smallskip
{\footnotesize\raggedright Notes. Modified policy minus LA-DRRC, using 40 matched episodes per load. Completion changes are percentage points.\par}

\end{table}

Universal diagnosis also admits several dispatch rules. Shortest diagnosis serves the eligible request with the shortest predicted service, while Recovery-$c\mu$ ranks requests by predicted recovery value per unit of service. Table~\ref{tab:alternativeordering} holds universal request and one-step recovery rules fixed. Reordering helps at high load, but all three rules still require every eligible failure to join the queue.

\begin{table}[htbp]
\centering
\caption{Universal Diagnosis with Alternative Queue Orderings\label{tab:alternativeordering}}
{\small\begin{tabular}{lrrrrrr}
\toprule
Policy & 0.50 & 0.80 & 1.05 & 1.25 & Complete (\%) & Expert time \\
\midrule
Diagnose-all FIFO & 0.554 & 0.444 & 0.290 & 0.181 & 37.19 & 0.897 \\
Shortest diagnosis & 0.555 & 0.474 & 0.383 & 0.328 & 55.15 & 0.881 \\
Recovery-$c\mu$ & 0.557 & 0.492 & 0.418 & 0.365 & 54.76 & 0.881 \\
\bottomrule
\end{tabular}\par}
\smallskip
{\footnotesize\raggedright Notes. The four load columns give net value. Completion and expert time refer to load 1.25. These are the same 40-episode tests used in Table~\ref{tab:closedloopmain}.\par}

\end{table}

\paragraph{Service quality, cost, and calibration.}
Panel A of Table~\ref{tab:sensitivityvalue} varies report accuracy and mean service at nominal load 1.05, adjusting arrival spacing to hold load fixed. The direct-cost rows use an independent set of paired seeds. At accuracy 0.55, gains over automated recovery are unresolved across all three service times. At accuracy 0.90, even slow service adds value. Increasing direct expert-time cost tenfold reduces the share of diagnosed offers from about 68\% to 2\%. Panel B varies one known condition at a time. The late-error kernel weights each stage by its order in the rebuilt suffix, while the persistent-error kernel retains the existing fault after a failed eligible recovery. Panel C changes the physical accuracy or duration without informing the controller. When actual accuracy is 0.55 and the rule assumes 0.78, using the available expert service lowers value relative to Automated DP.

\begin{table}[htbp]
\centering
\caption{Service Conditions and Errors in Controller Estimates\label{tab:sensitivityvalue}\label{tab:sensitivity}\label{tab:misspec}}
{\small\begin{tabular}{lrrr}
\toprule
Setting & LA value & vs. Auto [95\% CI] & vs. VOI [95\% CI] \\
\midrule
\multicolumn{4}{l}{\textit{A. Known report accuracy $a$, service multiplier $s$, and direct expert-time cost $c_H$}}\\
$a=0.55,\ s=0.75$ & 0.548 & 0.003 [-0.008, 0.013] & 0.024 [0.014, 0.034] \\
$a=0.55,\ s=1.0$ & 0.543 & -0.002 [-0.012, 0.007] & 0.024 [0.015, 0.033] \\
$a=0.55,\ s=1.5$ & 0.540 & -0.005 [-0.013, 0.002] & 0.033 [0.022, 0.043] \\
$a=0.78,\ s=0.75$ & 0.597 & 0.052 [0.039, 0.064] & 0.022 [0.011, 0.033] \\
$a=0.78,\ s=1.0$ & 0.583 & 0.038 [0.024, 0.052] & 0.014 [0.006, 0.022] \\
$a=0.78,\ s=1.5$ & 0.565 & 0.020 [0.008, 0.031] & 0.019 [0.010, 0.029] \\
$a=0.9,\ s=0.75$ & 0.624 & 0.079 [0.063, 0.095] & 0.017 [0.006, 0.027] \\
$a=0.9,\ s=1.0$ & 0.616 & 0.071 [0.059, 0.084] & 0.013 [0.004, 0.021] \\
$a=0.9,\ s=1.5$ & 0.598 & 0.053 [0.042, 0.064] & 0.018 [0.011, 0.025] \\
$c_H=0.015$ & 0.581 & 0.051 [0.039, 0.064] & 0.020 [0.010, 0.030] \\
$c_H=0.06$ & 0.555 & 0.025 [0.015, 0.036] & 0.023 [0.012, 0.034] \\
$c_H=0.15$ & 0.529 & -0.000 [-0.008, 0.007] & 0.028 [0.011, 0.044] \\
\midrule
\multicolumn{4}{l}{\textit{B. Known variability, deadlines, and error transitions}}\\
Service CV 0.25 & 0.588 & 0.045 [0.033, 0.057] & 0.026 [0.018, 0.035] \\
Service CV 0.75 & 0.580 & 0.037 [0.022, 0.051] & 0.026 [0.015, 0.036] \\
Deadline scale 0.80 & 0.546 & 0.048 [0.033, 0.062] & 0.017 [0.003, 0.030] \\
Deadline scale 1.20 & 0.601 & 0.041 [0.029, 0.053] & 0.016 [0.006, 0.025] \\
Late-error kernel & 0.597 & 0.046 [0.031, 0.060] & 0.023 [0.012, 0.034] \\
Persistent-error kernel & 0.587 & 0.037 [0.024, 0.051] & 0.026 [0.012, 0.040] \\
\midrule
\multicolumn{4}{l}{\textit{C. Physical and assumed settings differ (actual / assumed)}}\\
Accuracy 0.55 / 0.78 & 0.519 & -0.027 [-0.037, -0.017] & 0.016 [0.005, 0.028] \\
Accuracy 0.90 / 0.78 & 0.609 & 0.063 [0.052, 0.074] & 0.019 [0.009, 0.029] \\
Service 1.50 / 1.00 & 0.557 & 0.011 [0.000, 0.021] & 0.027 [0.016, 0.039] \\
\bottomrule
\end{tabular}\par}
\smallskip
{\footnotesize\raggedright Notes. Each row uses 20 paired episodes at nominal load 1.05. LA value is net value per offer; differences compare LA-DRRC with Automated DP and one-step VOI.\par}
\end{table}

\paragraph{Repeated investigations and concurrent operating changes.}
Longer later investigations may give the operator an incentive to obtain advice early. We keep initial arrival spacing fixed and inform the controller of the history-dependent service rule. Table~\ref{tab:history} shows that the gain over one-step VOI persists across all six cells, while the separate increment due to pricing later expert work ranges from $-0.0010$ to $0.0022$ and every interval covers zero.

\begin{table}[htbp]
\centering
\caption{History-Dependent Diagnostic Service\label{tab:history}}
{\scriptsize\begin{tabular}{lrrr}
\toprule
Setting & vs. Auto [95\% CI] & vs. VOI [95\% CI] & Price increment [95\% CI] \\
\midrule
$g=0.0,\ a=0.78$ & 0.057 [0.045, 0.069] & 0.029 [0.020, 0.039] & 0.000 [-0.003, 0.003] \\
$g=0.0,\ a=0.9$ & 0.075 [0.064, 0.085] & 0.022 [0.013, 0.032] & 0.001 [-0.002, 0.004] \\
$g=1.0,\ a=0.78$ & 0.053 [0.042, 0.063] & 0.036 [0.027, 0.045] & -0.001 [-0.002, 0.000] \\
$g=1.0,\ a=0.9$ & 0.068 [0.057, 0.079] & 0.020 [0.010, 0.029] & 0.002 [-0.001, 0.005] \\
$g=2.0,\ a=0.78$ & 0.052 [0.040, 0.063] & 0.035 [0.024, 0.046] & -0.001 [-0.002, 0.000] \\
$g=2.0,\ a=0.9$ & 0.068 [0.057, 0.080] & 0.017 [0.007, 0.027] & 0.000 [-0.001, 0.001] \\
\bottomrule
\end{tabular}\par}
\smallskip
{\footnotesize\raggedright Notes. Mean service is $s(x)(1+gf)$ after $f$ failed recoveries; $a$ is report accuracy. Price increment is LA-DRRC minus its no-future-work-price variant. Each cell uses 20 paired episodes.\par}

\end{table}

Table~\ref{tab:policyoodmain} changes several conditions together, giving the system one expert, more variable arrivals and service, and a less accurate report. At nominal loads 1.25 and 1.75, LA-DRRC's value differences from Automated DP are 0.003 and $-0.001$; both intervals cover zero. The share of offers diagnosed falls to 18.6\% and 6.4\%, respectively. These results trace the operating response as advice becomes harder to use.

\begin{table}[htbp]
\centering
\caption{Known Joint Change in Operating Conditions\label{tab:policyoodmain}}
{\small\begin{tabular}{lrrrrr}
\toprule
Load & Policy & Net value & Complete (\%) & Diagnosed (\%) & Expert time \\
\midrule
1.25 & Automated DP & 0.565 & 92.19 & 0.0 & 0.000 \\
1.25 & One-step VOI & 0.536 & 90.28 & 40.4 & 0.577 \\
1.25 & DRRC & 0.534 & 89.06 & 0.0 & 0.000 \\
1.25 & LA-DRRC & 0.568 & 92.20 & 18.6 & 0.257 \\
1.75 & Automated DP & 0.565 & 92.19 & 0.0 & 0.000 \\
1.75 & One-step VOI & 0.532 & 90.30 & 35.6 & 0.503 \\
1.75 & DRRC & 0.533 & 89.04 & 0.0 & 0.000 \\
1.75 & LA-DRRC & 0.564 & 92.37 & 6.4 & 0.080 \\
\bottomrule
\end{tabular}\par}
\smallskip
{\footnotesize\raggedright Notes. Thirty episodes per load, one expert, arrival CV 1.50, service multiplier 1.25, service CV 0.75, accuracy 0.60, and deadline multiplier 1.25. Estimators and coefficients are frozen; changed primitives are known.\par}

\end{table}

\FloatBarrier
\Needspace{4\baselineskip}
\subsection{Independent UCT and PPO benchmarks}
\label{app:repro-exact}\label{app:benchmarkdiagnostics}\label{app:ppo}\label{app:runtime}

UCT samples hidden errors from current beliefs, searches to depth eight, and uses six rollout steps under one-step VOI. It simulates fitted transitions among active workflows and aggregates observations by phase, focal workflow, report, attempts, rounded belief, feasible actions, and a 0.25-unit time bin. The exploration coefficient is 0.50; the leaf value sums positive one-recovery values across active workflows. Future arrivals occur in the evaluation system and are omitted from the planner. Budgets of 16, 64, and 256 use two development episodes of 24 offers per load at 0.80 and 1.25 (seeds 62900000--62900001). The corresponding development values are 0.368, 0.314 and 0.407. Testing uses the selected budget of 256 and planner seed equal to episode seed plus 777.

PPO trains on episodes with 180 offers, two experts, and loads sampled uniformly from the four test loads. The 49 observable features per candidate describe the action, phase, belief, task and queue conditions, and progress through the episode. All physically feasible actions enter the choice set, including those that DRRC would screen out at dispatch. A two-layer 64-unit tanh encoder feeds masked mean and maximum pooling and separate 64-unit actor and critic heads (28,098 parameters). PPO-Clip uses clipping 0.20, discount one, generalized-advantage parameter 0.95, 1,024 transitions per update, four epochs, minibatches of 128, and entropy and value-loss coefficients 0.01 and 0.50. The gradient norm is clipped at 0.50. The reward is the increment in operating value plus the change in an observable one-recovery potential, divided by five. The potential vanishes at termination, so shaping preserves the ordering of policies from a fixed initial state.

Learning rates $10^{-4}$ and $3\times 10^{-4}$ were compared on seeds 101 and 102, using 262,144 transitions per run and four validation episodes per load. Their validation values were 0.484 and 0.497. Five fresh seeds then trained for 524,288 transitions each at the selected higher rate; their selected updates were 384, 480, 288, 480 and 224. All methods were frozen before testing on the same 40 independent 180-offer episodes per load (seeds 83000000--83000039). The selected PPO policies made no diagnostic requests in their 800 test runs. In stochastic training traces, diagnosed offers fell from 5.38\% over the first 100 completed episodes per seed to 0.11\% over the last 100.

\begin{table}[htbp]
\centering
\caption{Matched UCT and PPO Comparison by Load\label{tab:ppoperformance}}
{\small\begin{tabular}{lrrrrrr}
\toprule
Policy & 0.50 & 0.80 & 1.05 & 1.25 & Complete (\%) & Expert time \\
\midrule
Automated DP & 0.541 & 0.541 & 0.541 & 0.541 & 87.18 & 0.000 \\
One-step VOI & 0.595 & 0.581 & 0.571 & 0.563 & 86.86 & 0.740 \\
UCT-256 & 0.465 & 0.458 & 0.446 & 0.439 & 78.49 & 0.234 \\
PPO & 0.521 & 0.521 & 0.522 & 0.522 & 86.56 & 0.000 \\
DRRC & 0.601 & 0.584 & 0.574 & 0.564 & 85.15 & 0.577 \\
LA-DRRC & 0.618 & 0.600 & 0.596 & 0.591 & 88.40 & 0.679 \\
\bottomrule
\end{tabular}\par}
\smallskip
{\footnotesize\raggedright Notes. Load columns report net value per offer. Completion and expert service per offer are at load 1.25. PPO averages over five training seeds. This sample is independent of Table~\ref{tab:closedloopmain}.\par}

\end{table}

Table~\ref{tab:ppoperformance} reports outcomes by load. For equally weighted averages, 10,000 bootstrap replicates (seed 86000000) resample the five training-seed indices and 40 matched episode indices independently, keeping each episode index together across loads and policies. Comparisons without PPO resample only episodes. LA-DRRC's four net-value gains over PPO are 0.097, 0.079, 0.075 and 0.070, with respective 95\% intervals $[0.087,0.107]$, $[0.069,0.089]$, $[0.065,0.085]$ and $[0.060,0.080]$. Section~\ref{sec:runtime} reports the pooled gains.

For decision-time comparisons, LA-DRRC and Diagnose-all FIFO each contribute ten states for every combination of E2, E3, E4 and loads 0.80 and 1.25. Every applicable policy sees the same copied states and cleared local caches. We time three calls per state and model, with randomized method order after warm-up on one AMD EPYC 9V74 logical core and one numerical thread. Table~\ref{tab:jointruntime} covers 3,480 calls including the five PPO networks. Timed work includes action enumeration, feature construction, inference and UCT's internal simulation. Loading, shared predictor inference, external state copying and offline training are excluded. Phase means enter the overall latency with equal weights over phases where the policy acts.

\begin{table}[htbp]
\centering
\caption{Shared-State Latency for All Benchmarks (Milliseconds)\label{tab:jointruntime}}
{\small\begin{tabular}{lrrrrr}
\toprule
Policy & E2 mean & E3 mean & E4 mean & E2/E4 median & E2/E4 95th \\
\midrule
Automated DP & 0.817 & -- & 0.641 & 0.661 & 0.981 \\
One-step VOI & 0.386 & 1.134 & 0.320 & 0.338 & 0.493 \\
UCT-256 & 701.592 & 701.957 & 675.269 & 549.997 & 2019.994 \\
PPO & 1.200 & 1.001 & 1.228 & 1.131 & 1.788 \\
DRRC & 0.619 & 1.158 & 0.491 & 0.509 & 0.860 \\
LA-DRRC & 0.613 & 1.109 & 0.481 & 0.480 & 0.845 \\
\bottomrule
\end{tabular}\par}
\smallskip
{\footnotesize\raggedright Notes. PPO pools five trained networks. Automated DP has no E3 decision.\par}

\end{table}

\FloatBarrier
\Needspace{4\baselineskip}
\subsection{Code repair across held-out repositories}
\label{app:empiricalppo}

\paragraph{Response profiles and evidence updates.}
Each training defect $j$ defines ten recorded binary outcomes $y_{ja}$ (five direct and five report-assisted routes) and a report category $\chi_j$. A regression tree fitted on training report text groups reports by their outcome vectors. For $n_{\rm tr}$ training defects, the smoothed outcomes and initial profile probabilities are
\[
\widetilde y_{ja}=\epsilon+(1-2\epsilon)y_{ja},\quad
\epsilon=(n_{\rm tr}+2)^{-1},\qquad
\pi_j^0=(1-\kappa)/n_{\rm tr}+\kappa\omega_{ij}.
\]
Here $\omega_{ij}$ averages normalized training-leaf membership across 96 trees fitted on the known problem description, faulty code, and test feedback to predict outcomes and report-category indicators. DRRC sets $\kappa=0$. With $K_{\rm rep}$ report categories, the likelihood is
\[
L_j(r)=(1-\rho_{\rm rep})/K_{\rm rep}
+\rho_{\rm rep}\mathbf1\{\chi_j=r\}.
\]
Failure of route $a$ multiplies $\pi_j$ by $1-\widetilde y_{ja}$; a returned report classified as $r$ multiplies it by $L_j(r)$. Normalization gives the posterior, and route $a$'s predicted success is $\sum_j\pi_j\widetilde y_{ja}$. Smoothing keeps failure likelihoods positive; smaller $\rho_{\rm rep}$ shrinks report information toward an uninformative likelihood. The contextual weight also blends predicted diagnostic work with its training mean: $\widehat s_i=\kappa\widehat s_i^{\rm ctx}+(1-\kappa)\bar s$. Both DRRC variants use common ridge estimates of route duration from input-size attributes. Fitting and parameter selection use training folds only.

\paragraph{Operating rules and work proxies.}
Unlike the queue commitment in the general model and controlled simulations, the software twin releases requests that fail the dispatch-time screen. All policies require enough predicted time for diagnosis and a recovery; One-step VOI, DRRC, and LA-DRRC additionally require positive diagnosis value at zero further wait, evaluated over the remaining one or two recovery attempts. A rejected request returns to a direct-only recovery-or-abandonment decision, with no immediate re-request. If that recovery fails, ordinary E2 decisions resume and may request diagnosis if no report has been obtained. No workflow repairs while queued, and each receives at most one diagnostic service. Diagnose-all FIFO and PPO use the physical screen without the additional positive-value requirement.

For a generated function of length $L$ characters, recovery takes $t=0.5+\max(250,L)/2000$ time units. Empty failed outputs use input length. A report of $L_H$ characters requires $s=0.5+L_H/1000$ units of expert time. Recovery costs $0.03t$, a started diagnosis costs $0.02+0.02s$, and waiting costs 0.005 per unit. Timely completion earns one unit. Arrivals are Poisson. Deadline slack is a Uniform$(3,7)$ draw times the recovery duration computed from input length. A started diagnosis occupies its server through service completion. Policies use work estimates fitted on training data.

Empirical PPO's 66 observable features exclude task and repository identifiers, reference repairs, pending reports, unchosen outcomes and future event times. Its action mask permits a request, direct repair or abandonment at E2, dispatch or idling at E3, and feasible repair or abandonment at E4. The same two-attempt, one-report limit applies to all policies.

PPO uses the synthetic benchmark's network, pooling, clipping, discount, generalized-advantage and loss settings, with 2,048 decisions per update across eight environments and minibatches of 256. Three repository-disjoint inner splits compare learning rates $10^{-4}$ and $3\times10^{-4}$ and four training budgets from 32,768 to 262,144 transitions. Each split refits the predictors and evaluates candidates on eight queues of 120 offers per load. All five outer folds select $10^{-4}$. Folds 0 and 3 select 131,072 transitions and the others select 65,536. Retraining on each outer training set with seeds 11, 23, 37, 53 and 71 yields 25 frozen policies. The matched queues in Section~\ref{sec:empirical} yield 480,000 PPO workflow records.

For the intervals in Section~\ref{sec:empirical}, we pair outcomes by held-out defect, weight loads equally and resample the 11 repository clusters. Comparisons with PPO use 5,000 bootstrap replicates that also resample training-seed indices, retaining each sampled index across folds. Fitted policies are held fixed.

\end{document}